\documentclass[11pt]{article}
\usepackage{amsfonts,amsmath,bm, xspace}
\usepackage{multirow,graphicx, algorithm,algorithmic,rotating}
\usepackage{url,cite}
\usepackage{fullpage}
\usepackage[small,bf]{caption}
\newtheorem{theorem}{Theorem}[section]

\newtheorem{proposition}[theorem]{Proposition}

\newcommand{\eat}[1]{}

\newcommand{\name}{\textsc{Hogwild!}\xspace}

\newcommand{\R}{\mathbb{R}}
\newcommand{\E}{\mathbb{E}}

\newcommand{\vct}[1]{\bm{#1}}
\newcommand{\mtx}[1]{\bm{#1}}

\newcommand{\norm}[1]{\left\lVert{#1}\right\rVert}

\newcommand{\fnorm}[1]{\norm{#1}_F}

\newcommand{\minimize}{\mbox{minimize}}
\newcommand{\st}{\mbox{subject to}}

\newcommand{\eq}[1]{(\ref{eq:#1})}

\numberwithin{equation}{section}

\title{\name: A Lock-Free Approach to Parallelizing Stochastic Gradient Descent}

\author{Feng Niu, Benjamin Recht, Christopher R\'{e} and Stephen J.~Wright\\
Computer Sciences Department, University of Wisconsin-Madison\\
1210 W Dayton St, Madison, WI 53706}

\date{June 2011}

\begin{document}

\maketitle

\vspace{-0.3in}

\begin{abstract}
Stochastic Gradient Descent (SGD) is a popular algorithm that can  achieve state-of-the-art performance on a variety of machine  learning tasks.  Several researchers have recently proposed schemes  to parallelize SGD, but all require performance-destroying memory locking and synchronization. This work aims to show using novel theoretical analysis, algorithms,  and implementation that SGD can be implemented {\em without any locking}.  We present an update scheme called \name which allows  processors access to shared memory with the possibility of  overwriting each other's work.  We show that when the  associated optimization problem is~\emph{sparse}, meaning most  gradient updates only modify small parts of the decision variable,  then \name achieves a nearly optimal rate of convergence.  We demonstrate experimentally that \name outperforms alternative schemes  that use locking by an order of magnitude.
\end{abstract}


\section{Introduction}


With its small memory footprint, robustness against noise, and rapid
learning rates, Stochastic Gradient Descent (SGD) has proved to be well suited to
data-intensive machine learning tasks~\cite{BertsekasNLP, Bottou08,
  Shalev-Shwartz08}. However, SGD's scalability is limited by its inherently sequential nature; it is difficult to parallelize. Nevertheless, the recent emergence of inexpensive multicore processors and mammoth, web-scale data sets has motivated researchers to develop several clever parallelization schemes for
SGD~\cite{BertsekasParallelBook,Langford09b,Dekel11,Zinkevich10,Duchi10}.
As many large data sets are currently pre-processed in a
MapReduce-like parallel-processing framework, much of the recent work
on parallel SGD has focused naturally on MapReduce implementations. MapReduce is
a powerful tool developed at Google for extracting information from
huge logs (e.g., ``find all the urls from a 100TB of Web data'') that was
designed to ensure fault tolerance and to simplify
the maintenance and programming of large clusters of
machines~\cite{Dean08}. But MapReduce is not ideally suited for
online, numerically intensive data analysis. Iterative computation is
difficult to express in MapReduce, and the overhead to ensure fault
tolerance can result in dismal throughput. Indeed, even Google
researchers themselves suggest that other systems, for example Dremel, are more
appropriate than MapReduce for data analysis tasks~\cite{Melnik10}.



For some data sets, the sheer size of the data dictates that one use a
cluster of machines. However, there are a host of
problems in which, after appropriate preprocessing, the data necessary
for statistical analysis may consist of a few terabytes or less. For
such problems, one can use a single inexpensive work station as
opposed to a hundred thousand dollar cluster. Multicore systems have
significant performance advantages, including  (1) low latency and high
throughput shared main memory (a processor in such a system can write
and read the shared physical memory at over 12GB/s with latency in the
tens of nanoseconds); and (2) high bandwidth off multiple disks (a
thousand-dollar RAID can pump data into main memory at over 1GB/s). In
contrast, a typical MapReduce setup will read incoming data at rates less than
tens of MB/s  due to frequent checkpointing for fault
tolerance. The high rates achievable by multicore systems move the bottlenecks
in parallel computation to synchronization (or locking) amongst the
processors~\cite{BerkeleyParallelTRShort,Fuller11}. Thus, to enable scalable data
analysis on a multicore machine, any performant solution must minimize
the overhead of locking.

In this work, we propose a simple strategy for eliminating the
overhead associated with locking: \emph{run SGD in parallel without
  locks}, a strategy that we call \name. In
\name, processors are allowed equal access to shared memory
and are able to update individual components of memory at will.
Such a lock-free scheme might appear doomed to fail as processors
could overwrite each other's progress.  However, when the data access
is~\emph{sparse}, meaning that individual SGD steps only modify a
small part of the decision variable, we show that memory overwrites are
rare and that they introduce barely any error into the computation when they do
occur.  We demonstrate both theoretically and experimentally a near
linear speedup with the number of processors on commonly occurring
sparse learning problems.

In Section~\ref{sec:graph}, we formalize a notion of
sparsity that is sufficient to guarantee such a speedup and provide
canonical examples of sparse machine learning problems 
in classification, collaborative filtering, and graph cuts.  Our notion of
sparsity allows us to provide theoretical guarantees of linear
speedups in Section~\ref{sec:fast-rates-theory}.  As a by-product of
our analysis, we also derive rates of convergence for algorithms with
constant stepsizes.  We demonstrate that robust $1/k$ convergence
rates are possible with constant stepsize schemes that implement an
exponential back-off in the constant over time.  This result is
interesting in of itself and shows that one need not settle for
$1/\sqrt{k}$ rates to ensure robustness in SGD algorithms.

In practice, we find that computational performance of a lock-free
procedure exceeds even our theoretical guarantees. We experimentally
compare lock-free SGD to several recently proposed methods.  We show
that all methods that propose memory locking are significantly slower
than their respective lock-free counterparts on a variety of machine learning
applications.

\section{Sparse Separable Cost Functions}\label{sec:graph}

\begin{figure}[t]
\centering
\includegraphics[width=3.75in]{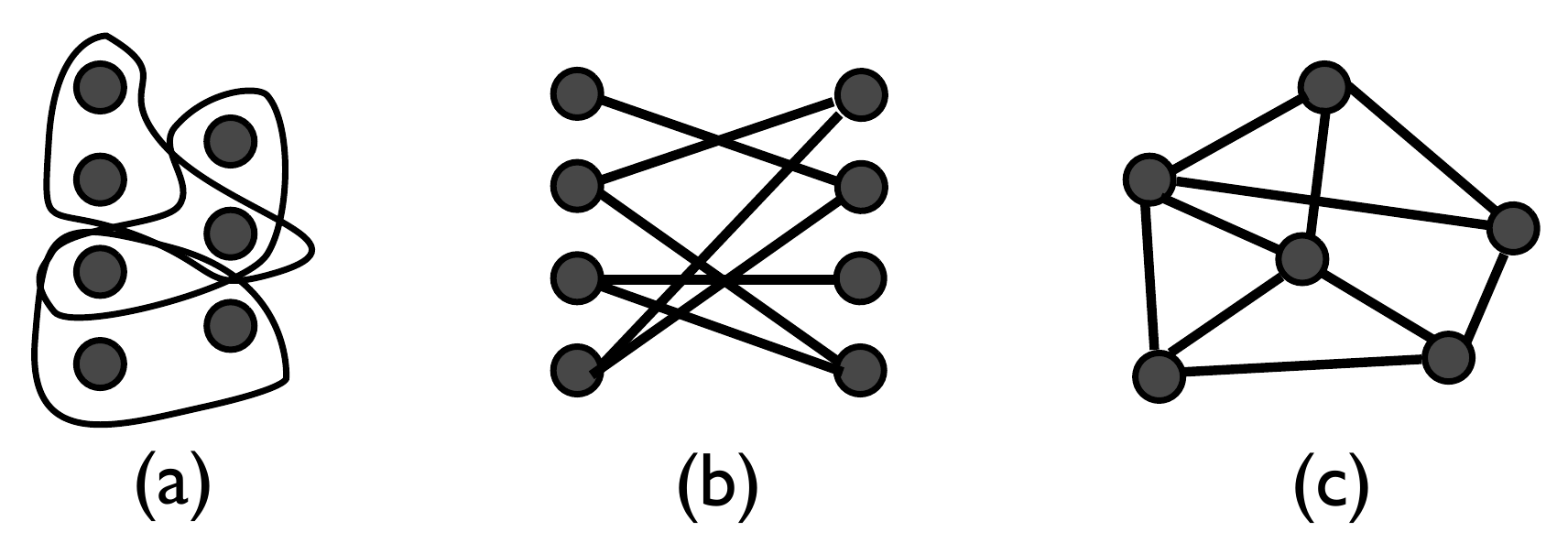}
\caption{\label{fig:graph-examples} Example graphs induced by cost function.  (a) A sparse SVM induces a hypergraph where each hyperedge corresponds to one example.  (b)  A matrix completion example induces a bipartite graph between the rows and columns with an edge between two nodes if an entry is revealed.  (c) The induced hypergraph in a graph-cut problem is simply the graph whose cuts we aim to find.}
\end{figure}

Our goal throughout is to minimize a function $f:X \subseteq \R^n\rightarrow \R$ of the form
\begin{equation}\label{eq:cost-function}
	f(x) = \sum_{e\in E} f_e (x_e)\,.
\end{equation}
Here $e$ denotes a small subset of $\{1,\ldots, n\}$ and $x_e$ denotes
the values of the vector $x$ on the coordinates indexed by $e$.  The
key observation that underlies our lock-free approach is that the natural cost functions 
associated with many machine learning problems of interest are
\emph{sparse} in the sense that $|E|$ and $n$ are both very large  but each individual $f_e$ acts only on a very small number of components of $x$.  That is, each subvector $x_e$
contains just a few components of $x$.

The cost function~\eq{cost-function} induces a \emph{hypergraph}
$G=(V,E)$ whose nodes are the individual components of $x$.  Each
subvector $x_e$ induces an edge in the graph $e\in E$ consisting of
some subset of nodes.  A few examples illustrate this concept.

\paragraph{Sparse SVM.} Suppose our goal is to fit a support vector machine to some data pairs $E=\{(z_1, y_1),\ldots,(z_{|E|}, y_{|E|})\}$ where $z\in \R^n$ and $y$ is a label for each $(z,y)\in E$.
\begin{equation}\label{eq:svm-cost}
	\minimize_x \sum_{\alpha \in E} \max(1 - y_\alpha x^T z_\alpha,0) +\lambda \|x\|_2^2\,,
\end{equation}
and we know~\emph{a priori} that the examples $z_\alpha$ are very sparse
(see for example \cite{Joachims06}).  To write this cost function in
the form of~\eq{cost-function}, let $e_\alpha$ denote the components which
are non-zero in $z_\alpha$ and let $d_u$ denote the number of training
examples which are non-zero in component $u$ ($u=1,2,\dotsc,n$).  Then
we can rewrite \eq{svm-cost} as
\begin{equation}\label{eq:svm-cost-sparse}
	\minimize_x \sum_{\alpha \in E} \left( \max(1 - y_\alpha x^T z_\alpha,0) +\lambda \sum_{u\in e_\alpha} \frac{x_u^2}{d_u} \right)\,.
\end{equation}
Each term in the sum~\eq{svm-cost-sparse}  depends only on the components of $x$ indexed by the set $e_\alpha$.

\paragraph{Matrix Completion.}  In the matrix completion problem, we
are provided entries of a low-rank, $n_r\times n_c$ matrix $\mtx{Z}$
from the index set $E$.  Such problems arise in collaborative
filtering, Euclidean distance estimation, and
clustering~\cite{Recht10,CandesRecht09,Lee10}. Our goal is to
reconstruct $\mtx{Z}$ from this sparse sampling of data.  A popular
heuristic recovers the estimate of $\mtx{Z}$ as a product $\mtx{L}\mtx{R}^*$ of factors
obtained from the following minimization:
\begin{equation}\label{eq:factor-nuclear}
	\minimize_{(\mtx{L},\mtx{R})} \,  \sum_{(u,v)\in E} (\vct{L}_{u} \vct{R}_v^*-Z_{uv})^2 + \tfrac{\mu}{2} \fnorm{\mtx{L}}^2 + \tfrac{\mu}{2} \fnorm{\mtx{R}}^2,
\end{equation}
where $\mtx{L}$ is $n_r\times r$, $\mtx{R}$ is $n_c\times r$ and
$\vct{L}_u$ (resp. $\vct{R}_v)$ denotes the $u$th (resp. $v$th) row of
$\mtx{L}$ (resp. $\mtx{R}$) ~\cite{mmmf,Recht10,Lee10}. To put
this problem in sparse form, i.e., as \eq{cost-function}, we write
\eq{factor-nuclear} as
\begin{equation*}\label{eq:factor-nuclear-collect}
	\minimize_{(\mtx{L},\mtx{R})} \sum_{(u,v)\in E } \left\{ (\vct{L}_{u} \vct{R}_v^*-Z_{uv})^2 +  \tfrac{\mu}{2|E_{u-}|} \|\vct{L}_u\|_F^2 + \tfrac{\mu}{2|E_{-v}|} \|\vct{R}_v\|_F^2\right\}
\end{equation*}
where $E_{u-} = \{v~:~(u,v)\in E\}$ and
$E_{-v}=\{u~:~(u,v)\in E\}$. 

\paragraph{Graph Cuts.}
Problems involving minimum cuts in graphs frequently arise in machine
learning (see~\cite{Boykov04} for a comprehensive survey).  In such
problems, we are given a sparse, nonnegative matrix $W$ which indexes similarity
between entities.  Our goal is to find a partition of the index set
$\{1,\ldots,n\}$ that best conforms to this similarity matrix.  Here
the graph structure is explicitly determined by the similarity matrix
$W$; arcs correspond to nonzero entries in $W$.  We want to match each
string to some list of $D$ entities.  Each node is associated with a
vector $x_i$ in the $D$-dimensional simplex $S_D= \{ \zeta \in \R^D~:~
\zeta_v\geq 0\,\,\, \sum_{v=1}^D \zeta_v = 1\}$.  Here, two-way cuts use $D=2$, but
multiway-cuts with tens of thousands of classes also arise in
entity resolution problems~\cite{Lee11}.  For example, we may have a
list of $n$ strings, and $W_{uv}$ might index the similarity of each string.  Several authors
(e.g.,~\cite{Cualinescu98}) propose to minimize the cost function
\begin{equation}\label{eq:graph-cut-cost}
	\minimize_{x}  \sum_{ (u,v) \in E } w_{uv} \|x_u - x_v\|_1 \quad \st \quad
x_v \in S_D\quad \mbox{for}~v=1,\ldots,n\,. 
\end{equation}


\vspace{.125in}

In all three of the preceding examples, the number of components involved in
a particular term $f_e$ is a small fraction of the total
number of entries.  We formalize this notion by defining the following
statistics of the hypergraph $G$:
\begin{align}
\label{eq:collide-prob}	\Omega := \max_{e\in E} |e|, ~~ \Delta  :=  \frac{\max_{1\leq v \leq n} |\{e\in E\,:\, v \in e\}|}{|E|}, ~~
	\rho := \frac{ \max_{e\in E}|\{\hat{e} \in E~: \hat{e}\cap e \neq \emptyset \}|}{|E|}\,. 
\end{align}
The quantity $\Omega$ simply quantifies the size of the hyper edges.
$\rho$ determines the maximum fraction of edges that intersect any
given edge. $\Delta$ determines the maximum fraction of edges that
intersect any variable.  $\rho$ is a measure of the sparsity of the
hypergraph, while $\Delta$ measures the node-regularity.  For our
examples, we can make the following observations about $\rho$ and
$\Delta$.
\begin{enumerate}
\item \textbf{Sparse SVM}.  $\Delta$ is simply the maximum frequency
  that any feature appears in an example, while $\rho$ measures how
  clustered the hypergraph is.  If some features are very common
  across the data set, then $\rho$ will be close to one.
	
\item \textbf{Matrix Completion}.  If we assume that the provided
  examples are sampled uniformly at random and we see more than $n_c
  \log(n_c)$ of them, then $\Delta \approx \tfrac{\log(n_r)}{n_r}$ and
  $\rho \approx \tfrac{2\log(n_r)}{n_r}$.  This follows from a
   \emph{coupon collector} argument~\cite{CandesRecht09}.
	
\item \textbf{Graph Cuts}.  $\Delta$ is the maximum degree divided by
  $|E|$, and $\rho$ is at most $2\Delta$.
\end{enumerate}

We now describe a simple protocol that achieves a linear speedup in the number of processors when $\Omega$, $\Delta$, and $\rho$ are relatively small.  

\section{The \name Algorithm}

Here we discuss the parallel processing setup.  We assume a shared
memory model with $p$ processors. The decision variable $x$ is
accessible to all processors. Each processor can read $x$, and can
contribute an update vector to $x$.  The vector $x$ is stored in
shared memory, and we assume that the componentwise addition operation
is atomic, that is
\[
	x_v \leftarrow x_v + a
\]
can be performed atomically by any processor for a scalar $a$ and
$v\in \{1,\ldots, n\}$.  This operation does not require a separate
locking structure on most modern hardware: such an operation is a
single atomic instruction on GPUs and DSPs, and it can be implemented via a
compare-and-exchange operation on a general purpose multicore
processor like the Intel Nehalem. In contrast, the operation of updating many components
at once requires an auxiliary locking structure.

Each processor then follows the procedure in Algorithm~\ref{alg:main-asynch}.   To fully describe the algorithm, let $b_v$ denote one of the standard basis elements in $\R^n$, with $v$ ranging from $1,\ldots, n$.  That is, $b_v$ is equal to $1$ on the $v$th component and $0$ otherwise.  Let $\mathcal{P}_v$ denote the Euclidean projection matrix onto the $v$th coordinate, i.e., $\mathcal{P}_v = b_v b_v^T$. $\mathcal{P}_v$ is a diagonal matrix equal to $1$ on the $v$th diagonal and zeros elsewhere.  Let $G_e(x)\in \R^n$ denote a gradient or subgradient of the function $f_e$ multiplied by $|E|$.   That is, we extend $f_e$ from a function on the coordinates of $e$ to all of $\R^n$ simply by ignoring the components in ${\neg e}$ (i.e., not in $e$).  Then 
\[
	|E|^{-1}G_e(x) \in \partial f_e (x).
\] 
Here, $G_e$ is equal to zero on the components in ${\neg e}$.   Using a sparse representation, we can calculate $G_e(x)$, only knowing the values of $x$ in the components indexed by $e$. Note that as a consequence of the uniform random sampling of $e$ from $E$, we have
\[
\E[G_e(x_e)] \in \partial f(x)\,.
\]

In Algorithm \ref{alg:psgd.worker}, each processor samples an term
$e\in E$ uniformly at random, computes the gradient of $f_e$ at $x_e$,
and then writes
\begin{equation} \label{update}
x_v \leftarrow x_v - \gamma  b_{v}^T G_{e}(x), \qquad \mbox{for each $v \in e$}
\end{equation}
Importantly,  note that the processor modifies only the variables indexed by $e$, leaving all of the components in ${\neg e}$ (i.e., not in $e$) alone.  We assume that the stepsize $\gamma$ is a fixed constant.  Even though the processors have no knowledge as to whether any of the other processors have modified $x$, we define $x_j$ to be the state of the decision variable $x$ after $j$ updates have been performed\footnote{Our notation overloads subscripts of $x$.  For clarity throughout, subscripts $i,j$, and $k$ refer to iteration counts, and $v$ and $e$ refer to components or subsets of components.}.
Since two processors can write to $x$ at the same time, we need to be
a bit careful with this definition, but we simply break ties at
random.  Note that $x_j$ is generally updated with a stale gradient, which is based on a value of $x$ read many
clock cycles earlier.  We use  $x_{k(j)}$ to denote the value of the
decision variable used to compute the gradient or subgradient that
yields the state $x_j$.

\begin{algorithm}[t]
  \caption{\name update for individual processors}
  \label{alg:psgd.worker}
  \begin{algorithmic}[1]
    \LOOP
    \STATE{Sample $e$ uniformly at random from $E$}
    \STATE{Read current state $x_e$ and evaluate $G_e (x)$}
    	\STATE{\textbf{for} $v\in e$ \textbf{do} $x_v \leftarrow x_v - \gamma b_{v}^T  G_{e}(x)$}
    \ENDLOOP
  \end{algorithmic}\label{alg:main-asynch}
\end{algorithm}

In what follows, we provide conditions under which this asynchronous, incremental gradient algorithm converges.  Moreover, we show that if the hypergraph induced by $f$ is isotropic and sparse, then this algorithm converges in nearly the same number of gradient steps as its serial counterpart.  Since we are running in parallel and without locks, this means that we get a nearly linear speedup in terms of the number of processors.

\section{Fast Rates for Lock-Free Parallelism}\label{sec:fast-rates-theory}

We now turn to our theoretical analysis of \name protocols.  To make the analysis tractable, we assume that we update with the following ``with replacement'' procedure:  each processor samples an edge $e$ uniformly at random and computes a subgradient of $f_e$ at the current value of the decision variable.  Then it chooses an $v\in e$ uniformly at random and updates
\begin{equation*}
	x_v \leftarrow x_v -\gamma |e| b_{v}^T G_e(x)
\end{equation*}
Note that the stepsize is a factor $|e|$ larger than the step in~(\ref{update}).  Also note that this update is completely equivalent to
\begin{equation}\label{eq:hw-wr-update}
	x \leftarrow x -\gamma |e| \mathcal{P}_{v}^T G_e(x)\,.
\end{equation}
This notation will be more convenient for the subsequent analysis.

This with replacement scheme assumes that a gradient is computed and then only one of its components is used to update the decision variable.  Such a scheme is computationally wasteful as the rest of the components of the gradient carry information for decreasing the cost.  Consequently, in practice and in our experiments, we perform a modification of this procedure.  We partition out the edges \emph{without replacement} to all of the processors at the beginning of each epoch.  The processors then perform \emph{full updates} of all of the components of each edge in their respective queues.  However, we emphasize again that we do not implement any locking mechanisms on any of the variables.  We do not analyze this ``without replacement'' procedure because no one has achieved tractable analyses for SGD in any without replacement sampling models. Indeed, to our knowledge, all analysis of without-replacement sampling yields rates that are comparable to a standard subgradient descent algorithm which takes steps along the full gradient of~\eq{cost-function}  (see, for example~\cite{Nedic00}).  That is, these analyses suggest that without-replacement sampling should require a factor of $|E|$ more steps than with-replacement sampling.  In practice, this worst case behavior is never observed.  In fact, it is conventional wisdom in machine learning that without-replacement sampling in stochastic gradient descent actually outperforms the with-replacement variants on which all of the analysis is based.

To state our theoretical results, we must describe several quantities that important in the analysis of our parallel stochastic gradient descent scheme.  We follow the notation and assumptions of Nemirovski \emph{et al}~\cite{Nemirovski09}. To simplify the analysis, we will assume that each $f_e$ in~\eq{cost-function} is a convex function.   
We assume Lipschitz continuous differentiability of $f$
with Lipschitz constant $L$:
\begin{equation} \label{eq:def.L}
\| \nabla f(x') - \nabla f(x) \| \le L \|x'-x \|, \;\; \forall \, x',
x \in X.  
\end{equation} 
We also assume $f$ is strongly convex with
modulus $c$. By this we mean that 
\begin{equation} \label{eq:def.c} 
f(x') \ge f(x) + (x'-x)^T\nabla f(x) + \frac{c}{2} \|x'-x\|^2, \;\;
\mbox{for all $x',x \in X$.}
\end{equation}
When $f$ is strongly convex, there exists a unique minimizer $x_\star$ and we denote $f_\star = f(x_\star)$.  We additionally assume that there exists a constant $M$ such that
\begin{equation} \label{eq:def.M} 
\| G_e(x_e) \|_2 \le M \;\; \mbox{almost surely for all $x \in X$}\,.
\end{equation}
We assume throughout that $\gamma c <1$.  (Indeed, when $\gamma c>1$, even the ordinary gradient descent algorithms will diverge.)

Our main results are summarized by the following
\begin{proposition}\label{prop:hw-main-result}
Suppose in Algorithm~\ref{alg:main-asynch} that the lag between when a gradient is computed and when it is used in step $j$ --- namely,  $j-k(j)$ --- is always less than or equal to $\tau$, and  $\gamma$ is defined to be
\begin{equation}\label{eq:hw-gamma-ub}
	\gamma = \frac{\vartheta \epsilon c}{2L M^2 \Omega \left( 1+6\rho \tau + 4\tau^2\Omega \Delta^{1/2}  \right)}\,.
\end{equation}
for some $\epsilon>0$ and $\vartheta\in(0,1)$.
Define $D_0:=\|x_0 - x_\star\|^2$ and let $k$ be an integer satisfying
\begin{equation}\label{eq:hw-global-iter-est}
	 	k \geq  \frac{2 L  M^2 \Omega \left( 1+ 6\tau \rho + 6 \tau^2 \Omega \Delta^{1/2}  \right) \log(LD_0/\epsilon)}{ c^2 \vartheta \epsilon }
\,.
\end{equation}
Then after $k$ component updates of $x$, we have $\E[ f(x_k) - f_\star ] \leq \epsilon$.
\end{proposition}

In the case that $\tau=0$, this reduces to precisely the rate achieved by the serial SGD protocol. A similar rate is achieved if $\tau=o(n^{1/4})$ as $\rho$ and $\Delta$ are typically both $o(1/n)$.  In our setting, $\tau$ is proportional to the number of processors, and hence as long as the number of processors is less $n^{1/4}$, we get nearly the same recursion as in the linear rate. 

We prove Proposition~\ref{prop:hw-main-result} in two steps in the Appendix.  First, we demonstrate that the sequence $a_j = \tfrac{1}{2} \E[ \|x_j - x_\star\|^2]$ satisfies a recursion of the form $a_j \leq (1-c_r\gamma) (a_{j+1}-a_\infty) + a_\infty$ for some constant $a_\infty$ that depends on many of the algorithm parameters but not on the state, and some constant $c_r < c$. This $c_r$ is an ``effective curvature'' for the problem which is smaller that the true curvature $c$ because of the errors introduced by our update rule.   Using the fact that $c_r\gamma <1$, we will show in Section~\ref{sec:css-rates} how to determine an upper bound on $k$ for which $a_k \leq \epsilon/L$.  Proposition~\ref{prop:hw-main-result} then follows because $E[f(x_k)-f(x_\star)] \leq L a_k$ since the gradient of $f$ is Lipschitz.  A full proof is provided in the appendix.

Note that up to the $\log(1/\epsilon)$ term in~\eq{hw-global-iter-est}, our analysis nearly provides a $1/k$ rate of convergence for a constant stepsize SGD scheme, both in the serial and parallel cases.   Moreover, note that our rate of convergence is fairly robust to error in the value of $c$; we pay linearly for our underestimate of the curvature of $f$.  In contrast,  Nemirovski \emph{et al} demonstrate that when the stepsize is inversely proportional to the iteration counter, an overestimate of $c$ can result in exponential slow-down~\cite{Nemirovski09}!  We now turn to demonstrating that we can eliminate the log term from~\eq{hw-global-iter-est} by a slightly more complicated protocol where the stepsize is slowly decreased after a large number of iterations.

\section{Robust $1/k$ rates.}\label{sec:css-rates} 

Suppose we run Algorithm~\ref{alg:main-asynch} for a fixed number of gradient updates $K$ with stepsize $\gamma<1/c$.  Then, we wait for the threads to coalesce, reduce $\gamma$ by a constant factor $\beta \in (0,1)$, and run  for $\beta^{-1} K$ iterations.  In some sense, this piecewise constant stepsize protocol approximates a $1/k$ diminishing stepsize.  The main difference with the following analysis from previous work is that our stepsizes are always less than $1/c$ in contrast to beginning with very large stepsizes.  Always working with small stepsizes allows us to avoid the possible exponential slow-downs that occur with standard diminishing stepsize schemes.

To be precise, suppose $a_k$ is any sequence of real numbers satisfying
\begin{equation}\label{eq:hw-base-recursion}
	a_{k+1} \leq (1-c_r \gamma) (a_k - a_\infty(\gamma)) + a_\infty(\gamma)
\end{equation}
where $a_\infty$ is some non-negative function of $\gamma$ satisfying 
\[
	a_\infty(\gamma) \leq \gamma B
\]
and $c_r$ and $B$ are constants. This recursion underlies many convergence proofs for SGD where $a_k$ denotes the distance to the optimal solution after $k$ iterations. We will derive appropriate constants for \name in the Appendix.  We will also discuss below what these constants are for standard stochastic gradient descent algorithms.

Factoring out the dependence on $\gamma$ will be useful in what follows. Unwrapping~\eq{hw-base-recursion} we have
\[
	a_{k} \leq (1-c_r\gamma)^k (a_0 - a_\infty(\gamma)) +a_\infty(\gamma) \,.
\]
Suppose we want this quantity to be less than $\epsilon$.  It is sufficient that both terms are less than $\epsilon/2$.  For the second term, this means that it is sufficient to set
\begin{equation}\label{eq:hw-ball-constraint1}
	\gamma \leq \frac{\epsilon}{2B}\,.
\end{equation}
For the first term, we then need
\[
	\left(1-\gamma c_r \right)^k a_0 \leq \epsilon/2
\]
which holds if
\begin{equation}\label{eq:hw-one-epoch}
	 k \geq \frac{\log(2a_0/\epsilon)}{\gamma c_r}\,.
\end{equation}

By~\eq{hw-ball-constraint1}, we should pick $\gamma = \frac{\epsilon\vartheta}{2B}$ for $\vartheta \in (0,1]$.  Combining this with~\eq{hw-one-epoch} tells us that after
\[
	 k \geq \frac{2B\log(2a_0/\epsilon)}{\vartheta \epsilon c_r}
\]
iterations we will have $a_k \leq \epsilon$.  This right off the bat almost gives us a $1/k$ rate, modulo the $\log(1/\epsilon)$ factor.  

To eliminate the log factor, we can implement a backoff scheme where we reduce the stepsize by a constant factor after several iterations.   This backoff scheme will have two phases: the first phase will consist of converging to the ball about $x_\star$ of squared radius less than $\frac{2B}{c_r}$ at an exponential rate.  Then we will converge to $x_\star$ by shrinking the stepsize.

To calculate the number of iterates required to get inside a ball of squared radius $\frac{2B}{c_r}$, suppose the initial stepsize is chosen as $\gamma= \frac{\vartheta}{c_r}$ ($0<\vartheta<1$). This choice of stepsize guarantees that the $a_k$ converge to $a_\infty$.  We use the parameter $\vartheta$ to demonstrate that we do not suffer much for underestimating the optimal stepsize (i.e., $\vartheta=1$) in our algorithms.  Using~\eq{hw-one-epoch} we find that
\begin{equation}\label{eq:hw-linear-phase}
	 k \geq  \vartheta^{-1} \log\left(\frac{a_0 c_r}{\vartheta B}\right)
\end{equation}
iterations are sufficient to converge to this ball.  Note that this is a linear rate of convergence.

Now assume that $a_0 < \frac{2\vartheta B}{c_r}$.  Let's reduce the stepsize by a factor of $\beta$ each epoch.  This reduces the achieved $\epsilon$ by a factor of $\beta$.  Thus, after $\log_\beta(a_0/\epsilon)$ epochs, we will be at accuracy $\epsilon$.  The total number of iterations required is then the sum of terms with the form~\eq{hw-one-epoch}, with $a_0$ set to be the radius achieved by the previous epoch and $\epsilon$ set to be $\beta$ times this $a_0$.  Hence, for epoch number $\nu$, the initial distance is $\beta^{\nu-1}a_0$ and the final radius is $\beta^\nu$. Summing over all of the epochs (except for the initial phase) gives
\begin{equation}\label{eq:hw-poly-phase}
\begin{aligned}
	\sum_{k=1}^{\log_\beta(a_0/\epsilon)} \frac{\log(2/\beta)}{\vartheta \beta^k} &=
	\frac{\log(2/\beta)}{\vartheta} \sum_{k=1}^{\log_\beta(a_0/\epsilon)} \beta^{-k}\\
&=	\frac{\log(2/\beta)}{\vartheta}\frac{\beta^{-1}(a_0/\epsilon) - 1}{\beta^{-1}-1}\\
&\leq\frac{a_0}{\vartheta\epsilon}\frac{\log(2/\beta)}{1-\beta}\\ 
&\leq	\frac{2B}{ c_r \epsilon}\frac{\log(2/\beta)}{1-\beta}\,.
\end{aligned}
\end{equation}
This expression is minimized by selecting a backoff parameter $\approx 0.37$.  Also, note that when we reduce the stepsize by $\beta$, we need to run for $\beta^{-1}$ more iterations. 

Combining~\eq{hw-linear-phase} and~\eq{hw-poly-phase}, we estimate a total number of iterations equal to
\[
	k \geq \vartheta^{-1} \log\left(\frac{a_0 c_r}{\vartheta B}\right)+\frac{2B}{c_r \epsilon}\frac{\log(2/\beta)}{1-\beta}
\]
are sufficient to guarantee that $a_k \leq \epsilon$.

Rearranging terms, the following two expressions give $\epsilon$ in terms of all of the algorithm parameters:
 \begin{equation}\label{eq:css-error-estimate}
	\epsilon \leq \frac{2\log(2/\beta)}{1-\beta}\cdot \frac{B}{c_r} \cdot\frac{1}{k-\vartheta^{-1} \log\left(\frac{a_0c_r}{\vartheta B}\right)}\,.
\end{equation}

\subsection{Consequences for serial SGD}
Let us compare the results of this constant step-size protocol to one where the stepsize at iteration $k$ is set to be $\gamma_0/k$ for some initial step size $\gamma$ for the standard (serial) incremental gradient algorithm applied to~\eq{cost-function}.   Nemirovski \emph{et al} ~\cite{Nemirovski09} show that the expected squared distance to the optimal solution, $a_k$, satisfies
\[
	a_{k+1} \leq (1-2c\gamma_k) a_k + \tfrac{1}{2}\gamma_k^2 M^2\,.
\]
We can put this recursion in the form~\eq{hw-base-recursion} by setting $\gamma_k = \gamma$, $c_r = 2c$, $B = \tfrac{M^2}{4c}$, and $a_\infty =\tfrac{\gamma M^2}{4c}$.

The authors of~\cite{Nemirovski09} demonstrate that a large step size: $\gamma_k = \tfrac{\Theta}{2ck}$ with $\Theta>1$  yields a bound
\[
	a_k \leq\frac{1}{k} \max\left\{ \frac{M^2}{c^2} \cdot \frac{\Theta^2}{4\Theta-4},  D_0    \right\}
\]

On the other hand, a constant step size protocol achieves
 \[
	a_k \leq \frac{\log(2/\beta)}{4(1-\beta)}\cdot \frac{M^2}{c^2} \cdot\frac{1}{k-\vartheta^{-1} \log\left(\frac{4D_0c^2}{\vartheta M^2}\right)}\,.
\]
This bound is obtained by plugging the algorithm parameters into~\eq{css-error-estimate} and letting $D_0 = 2 a_0$.

Note that both bounds have asymptotically the same dependence on $M$, $c$, and $k$.  The expression
\[
	\frac{\log(2/\beta)}{4(1-\beta)}
\]
is minimized when $\beta \approx 0.37$ and is equal to $1.34$.
The expression
\[
\frac{\Theta^2}{4\Theta-4}
\]
is minimized when $\Theta = 2$ and is equal to $1$ at this minimum.  So the leading constant is slightly worse in the constant stepsize protocol when all of the parameters are set optimally. However, if $D_0 \geq M^2/c^2$,  the $1/k$ protocol has error proportional to $D_0$, but our constant stepsize protocol still has only a logarithmic dependence on the initial distance. Moreover, the constant stepsize scheme is much more robust to overestimates of the curvature parameter $c$.  For the $1/k$ protocols, if one overestimates the curvature (corresponding to a small value of $\Theta$), one can get arbitrarily slow rates of convergence.  An simple, one dimensional example in~\cite{Nemirovski09} shows that $\Theta=0.2$ can yield a convergence rate of $k^{-1/5}$.  In our scheme, $\vartheta=0.2$ simply increases the number of iterations by a factor of $5$.

The proposed fix in~\cite{Nemirovski09} for the sensitivity to curvature estimates results in asymptotically slower convergence rates of $1/\sqrt{k}$.  It is important to note that we  need not settle for these slower rates and can still achieve robust convergence at $1/k$ rates.

\subsection{Parallel Implementation of a Backoff Scheme}
The scheme described about results in a $1/k$ rate of convergence for \name with the only synchronization overhead occurring at the end of each ``round'' or ``epoch'' of iteration.   When implementing a backoff scheme for \name, the processors have to agree on when to reduce the stepsize.  One simple scheme for this is to run all of the processors for a fixed number of iterations, wait for all of the threads to complete, and then globally reduce the stepsize in a master thread.  We note that one can eliminate the need for the threads to coalesce by sending out-of-band messages to the processors to signal when to reduce $\gamma$.  This complicates the theoretical analysis as there may be times when different processors are running with different
stepsizes, but in practice could allow one to avoid synchronization costs.  We do not implement this scheme, and so do not analyze this idea further.

\section{Related Work}

Most schemes for parallelizing stochastic gradient descent are
variants of ideas presented in the seminal text 
by Bertsekas and Tsitsiklis~\cite{BertsekasParallelBook}.  For
instance, in this text, they describe using stale gradient updates
computed across many computers in a master-worker setting and
describe settings where different processors control access to
particular components of the decision variable.  They prove global
convergence of these approaches, but do not provide rates of
convergence (This is one way in which our work extends this prior
research).  These authors also show that SGD convergence is robust to a
variety of models of delay in computation and communication
in~\cite{Tsitsiklis86}.

We also note that constant stepsize protocols with backoff procedures are canonical in SGD practice, but perhaps not in theory.  Some theoretical work which has at least demonstrated convergence of these protocols can be found in~\cite{Luo94, Tseng98}.  These works do not establish the $1/k$ rates which we provided above.

Recently, a variety of parallel schemes have been proposed in a variety of contexts.  In MapReduce settings, Zinkevich et al proposed running many instances of stochastic gradient descent on different machines and averaging their output~\cite{Zinkevich10}.  Though the authors claim this method can reduce both the variance of their estimate and the overall bias, we show in our experiments that for the sorts of problems we are concerned with, this method does not outperform a serial scheme.

Schemes involving the averaging of gradients via a distributed protocol have also been proposed by several authors~\cite{Dekel11,Duchi10}.  While these methods do achieve linear speedups, they are difficult to implement efficiently on multicore machines as they require massive communication overhead.  Distributed averaging of gradients requires message passing between the cores, and the cores need to synchronize frequently in order to compute reasonable gradient averages.

The work most closely related to our own is a round-robin scheme
proposed by Langford et al~\cite{Langford09b}.  In this scheme, the
processors are ordered and each update the decision variable in order.
When the time required to lock memory for writing is dwarfed by the
gradient computation time, this method results in a linear speedup, as
the errors induced by the lag in the gradients are not too severe.
However, we note that in many applications of interest in machine
learning, gradient computation time is incredibly fast, and we now
demonstrate that in a variety of applications, \name outperforms such
a round-robin approach by an order of magnitude.

\section{Experiments}

\begin{figure}
\small
\begin{center}
\begin{tabular}{|c|c|rrr|rrr|rrr|}
\multicolumn{5}{c}{ }&\multicolumn{3}{c}{\name} &\multicolumn{3}{c}{\textsc{Round Robin}}\\
\hline
\multirow{2}{*}{type}&data & \multicolumn{1}{c}{size} &  \multicolumn{1}{c}{$\rho$} &  \multicolumn{1}{c|}{$\Delta$} & \multicolumn{1}{c}{time} & \multicolumn{1}{c}{train} &  \multicolumn{1}{c|}{test}  & \multicolumn{1}{c}{time} &  \multicolumn{1}{c}{train} &  \multicolumn{1}{c|}{test} \\
 &set & \multicolumn{1}{c}{(GB)} & & & \multicolumn{1}{c}{(s)} &  \multicolumn{1}{c}{error} &  \multicolumn{1}{c|}{error} & \multicolumn{1}{c}{(s)}    &  \multicolumn{1}{c}{error} &  \multicolumn{1}{c|}{error} \\
\hline
\textbf{SVM}
& RCV1 & 0.9 &0.44 & 1.0 & 9.5 & 0.297 & 0.339 & 61.8 & 0.297 & 0.339 \\
\hline
\multirow{3}{*}{\textbf{MC}}
& Netflix   & 1.5  & 2.5e-3 & 2.3e-3 & 301.0 & 0.754 & 0.928 & 2569.1 &  0.754 & 0.927\\
& KDD  & 3.9  & 3.0e-3 & 1.8e-3 & 877.5 & 19.5 & 22.6 &  7139.0 & 19.5 &  22.6\\
& Jumbo   & 30 & 2.6e-7 & 1.4e-7 &9453.5 & 0.031 & 0.013 & \textbf{N/A} &  \textbf{N/A} &  \textbf{N/A}\\
\hline
\multirow{2}{*}{\textbf{Cuts}} 
& DBLife  & 3e-3  & 8.6e-3 & 4.3e-3 & 230.0 & 10.6 &  \textbf{N/A} & 413.5 & 10.5 &  \textbf{N/A}\\
& Abdomen  & 18  & 9.2e-4 & 9.2e-4 & 1181.4 & 3.99 &  \textbf{N/A} & 7467.25 & 3.99 &  \textbf{N/A}\\
\hline
\end{tabular}
\caption{Comparison of wall clock time across of \name and RR. Each algorithm is run for $20$ epochs and parallelized over 10 cores.}\label{fig:table-compare}
\end{center}
\end{figure}

We ran numerical experiments on a variety of machine learning tasks, and compared against a round-robin approach proposed in~\cite{Langford09b} and implemented in Vowpal Wabbit~\cite{VowpalWabbit}.  We refer to this approach as RR. To be as fair as possible to prior art, we hand coded RR to be nearly identical to the \name approach, with the only difference being the schedule for how the gradients are updated.  One notable change in RR from the Vowpal Wabbit software release is
that we optimized RR's locking and signaling mechanisms to
use spinlocks and busy waits (there is no need for generic signaling
to implement round robin). We verified that this optimization results
in nearly an order of magnitude increase in wall clock time for all
problems that we discuss.

We also compare against a model which we call AIG which can be seen as
a middle ground between RR and \name.  AIG runs a protocol identical
to \name except that it locks all of the variables in $e$ in before
and after the \textbf{for} loop on line 4 of
Algorithm~\ref{alg:psgd.worker}.  Our experiments demonstrate that
even this fine-grained locking induces undesirable slow-downs.

All of the experiments were coded in C++ are run on an identical configuration: a dual Xeon X650 CPUs (6 cores each x 2 hyperthreading) machine with 24GB of RAM and a software RAID-0 over 7 2TB Seagate Constellation 7200RPM disks.  The kernel is Linux 2.6.18-128.  We never use more than 2GB of memory. All training data is stored on a seven-disk raid 0. We implemented a custom file scanner to demonstrate the speed of reading data sets of disk into small shared memory.  This allows us to read data from the raid at a rate of nearly 1GB/s.

All of the experiments use a constant stepsize $\gamma$ which is diminished by a factor $\beta$ at the end of each pass over the training set.  We run all experiments for 20 such passes, even though less epochs are often sufficient for convergence.   We show results for the largest value of the learning rate $\gamma$ which converges and we use $\beta=0.9$ throughout.  We note that the results look the same across a large range of $(\gamma,\beta)$ pairs and that all three parallelization schemes achieve train and test errors within a few percent of one another.  We present experiments on the classes of problems described in Section~\ref{sec:graph}.

\begin{figure}
\centering
\begin{tabular}{ccc}
\includegraphics[width=50mm]{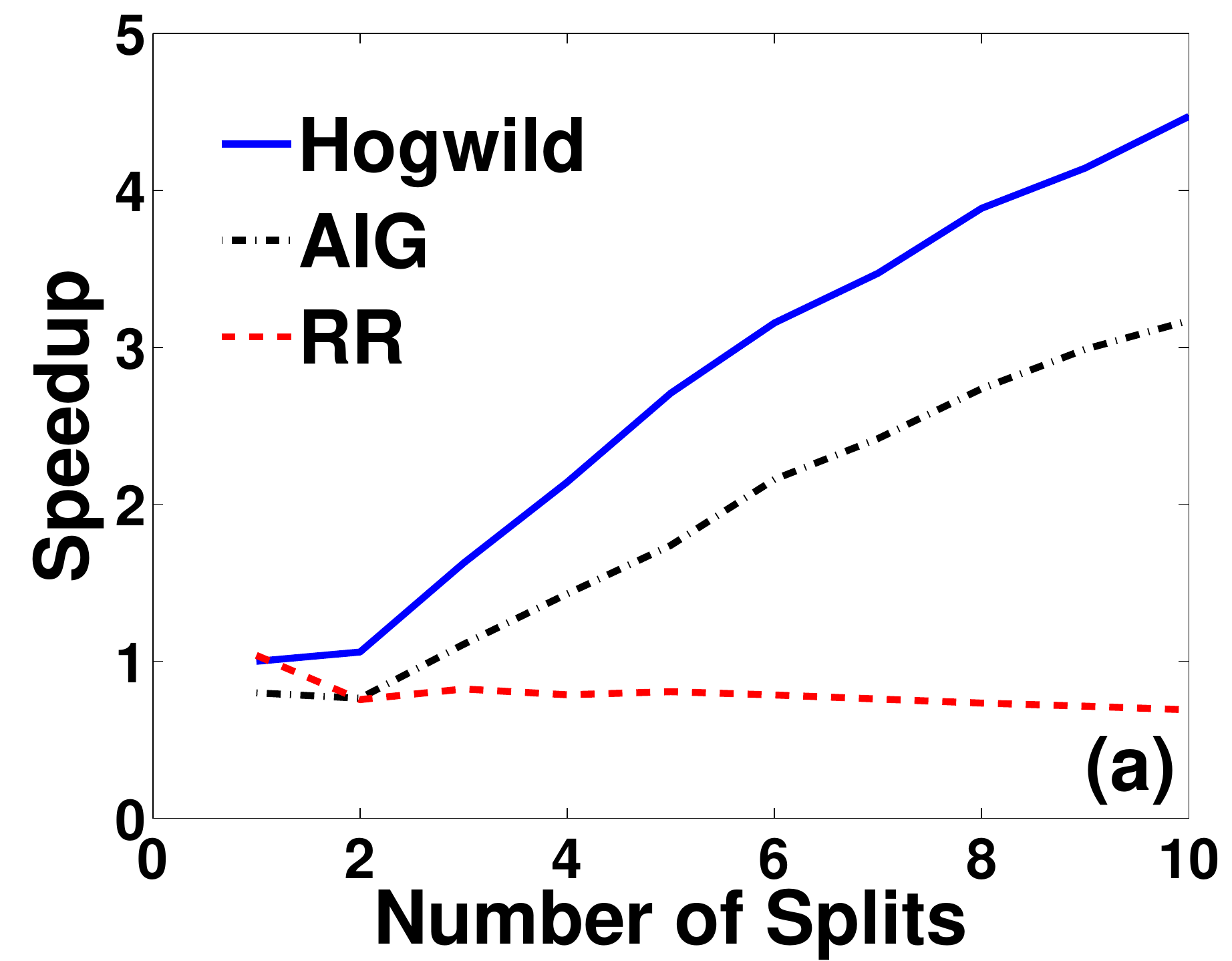} & 
\includegraphics[width=50mm]{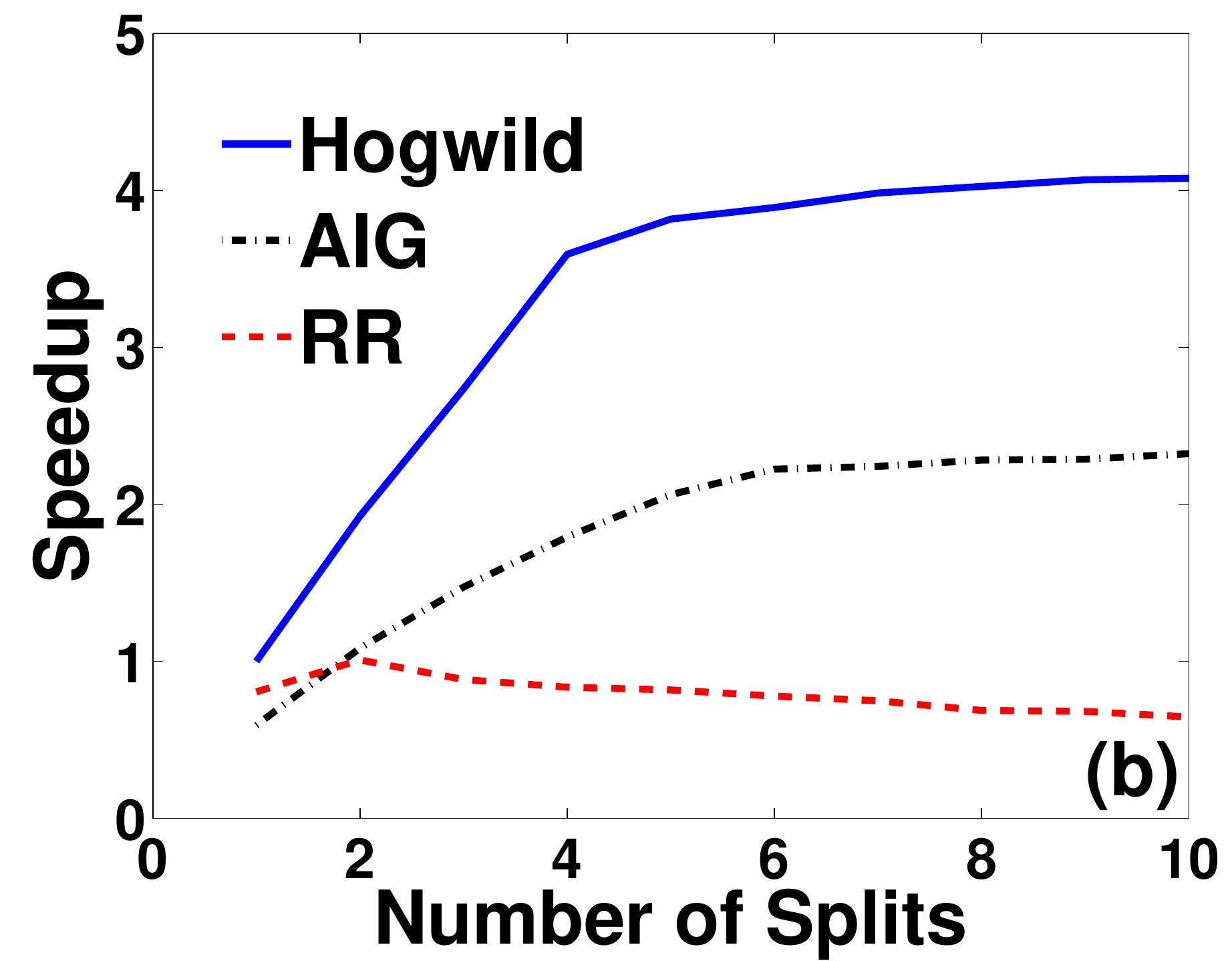} & 
\includegraphics[width=50mm]{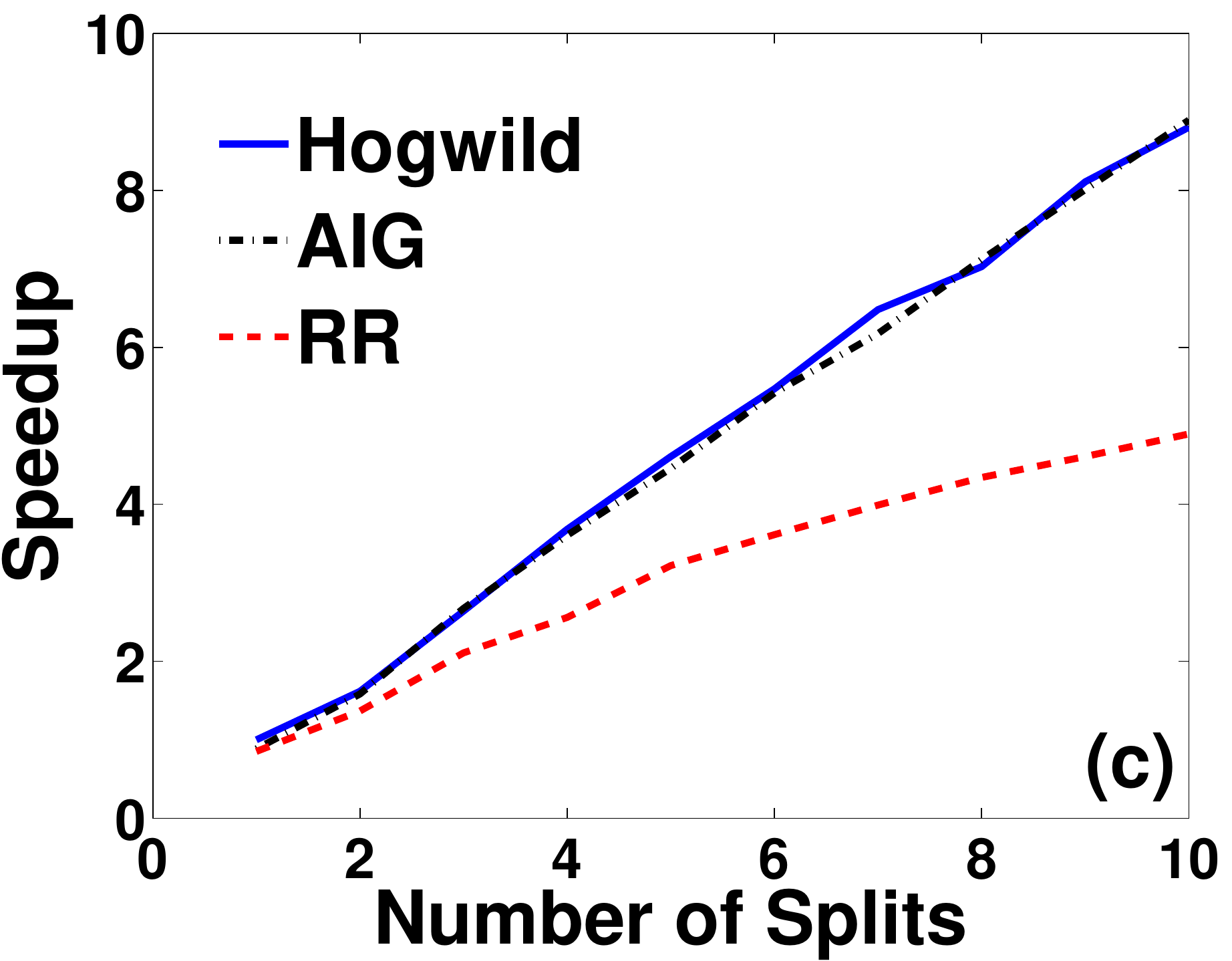} 
\end{tabular}
\caption{Total CPU time versus number of threads for (a) RCV1, (b) Abdomen, and (c) DBLife.
\label{fig:speedups}}
\end{figure}

\textbf{Sparse SVM.}  We tested our sparse SVM implementation on the Reuters RCV1 data set on the binary
text classification task CCAT~\cite{Lewis04}. There are 804,414
examples split into 23,149 training and 781,265 test examples, and
there are 47,236 features. We swapped the training set and the test set for our experiments to
demonstrate the scalability of the parallel multicore algorithms. In
this example, $\rho=0.44$ and $\Delta=1.0$---large
values that suggest  a bad case for \name.  Nevertheless, in Figure~\ref{fig:speedups}(a),
we see that \name is able to achieve a factor of 3 speedup with while
RR gets worse as more threads are added.  Indeed, for fast gradients,
RR is worse than a serial implementation.

\begin{figure}
\centering
\begin{tabular}{ccc}
\includegraphics[width=50mm]{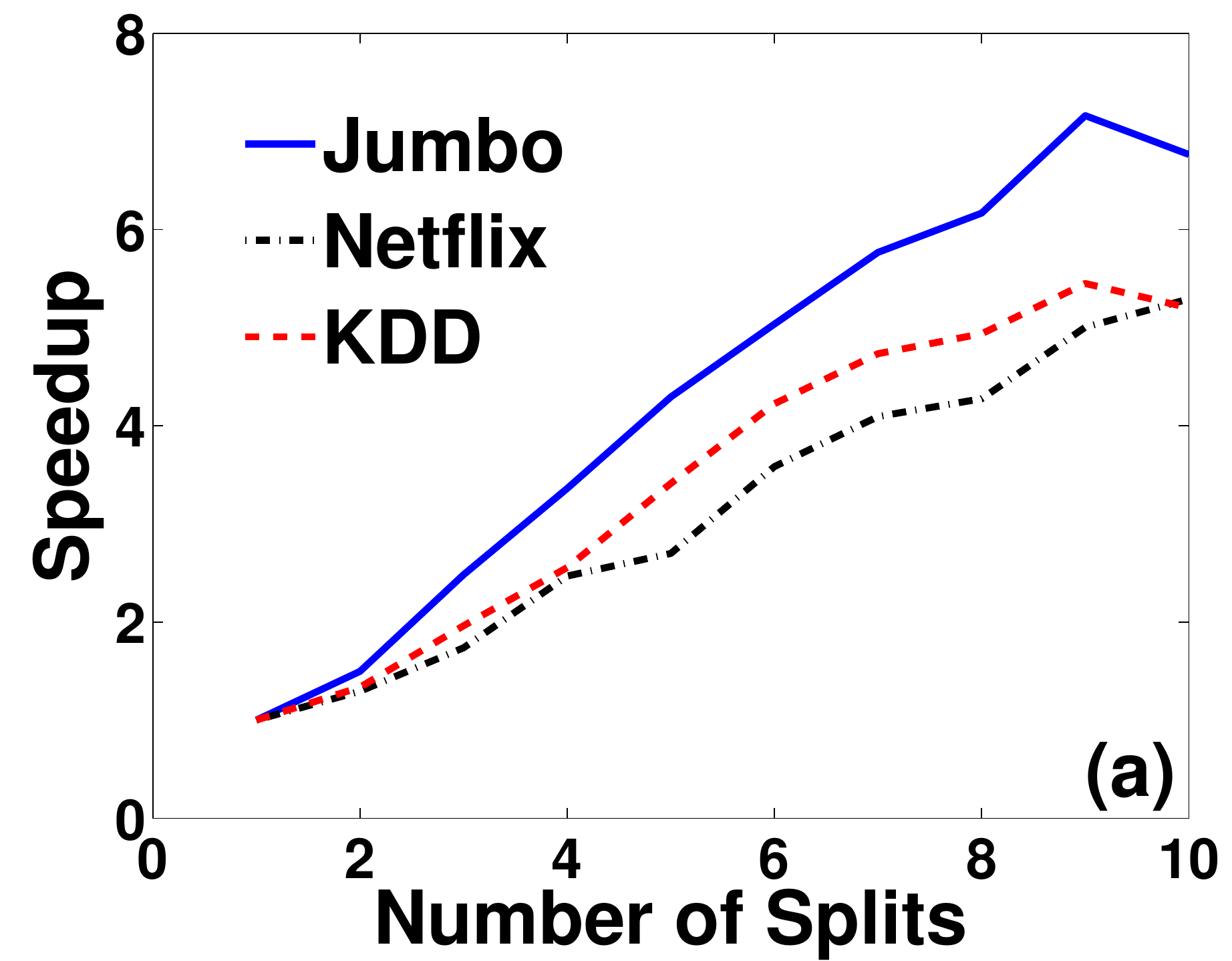} & 
\includegraphics[width=50mm]{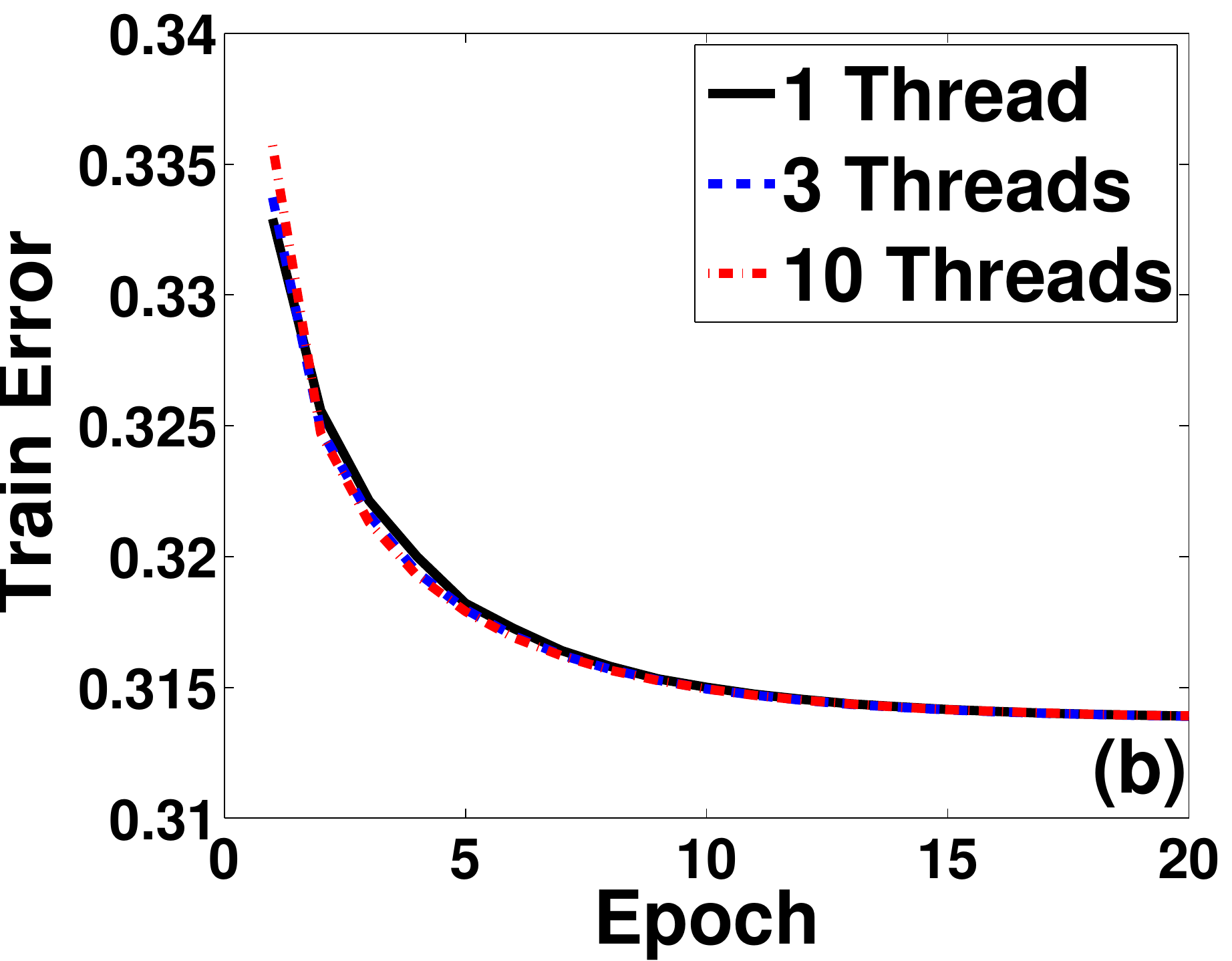} & 
\includegraphics[width=50mm]{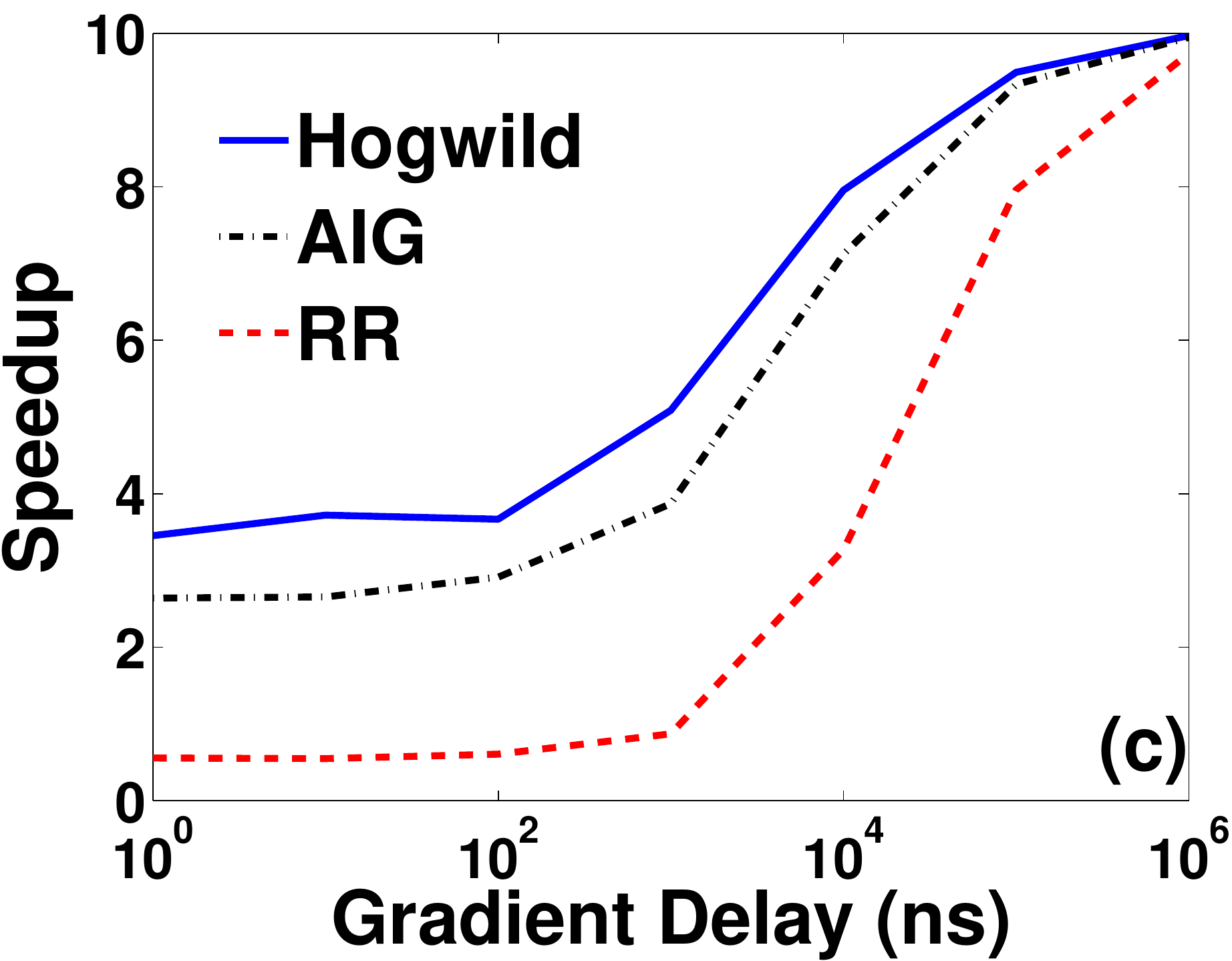} 
\end{tabular}
\caption{(a) Speedup for the three matrix completion problems with \name.  In all three cases, massive speedup is achieved via parallelism.  (b) The training error at the end of each epoch of SVM training on RCV1 for the averaging algorithm~\cite{Zinkevich10}.  (c) Speedup achieved over serial method for various  levels of delays (measured in nanoseconds).
\label{fig:misc}}
\end{figure}

For this data set, we also implemented the approach
in~\cite{Zinkevich10} which runs multiple SGD runs in parallel and
averages their output.  In Figure~\ref{fig:misc}(b), we display at the
train error of the ensemble average across parallel threads at the end
of each pass over the data.   We note that the threads only communicate at the very
end of the computation, but we want to demonstrate the effect of parallelization on train
error.   Each of the parallel threads touches every data example in each pass. Thus, the $10$ thread run does $10$x more gradient computations than
the serial version. Here,  the error is the same whether we run in serial
or with ten instances.  We conclude that on this problem, there is no advantage to
running in parallel with this averaging scheme.

\textbf{Matrix Completion.}   We ran \name on three very large matrix completion problems.  The Netflix Prize data set has 17,770 rows, 480,189 columns, and 100,198,805 revealed entries.  The KDD Cup 2011 (task 2) data set has 624,961 rows, 1,000,990, columns and 252,800,275 revealed entries.   We also synthesized  a low-rank matrix with rank $10$, 1e7 rows and columns, and 2e9 revealed entries.   We refer to this instance as ``Jumbo.''  In this synthetic example, $\rho$ and $\Delta$ are both around 1e-7.  These values  contrast sharply with the real data sets where $\rho$ and $\Delta$ are both on the order of 1e-3.

\begin{figure}
\centering
\begin{tabular}{ccc}
\includegraphics[width=50mm]{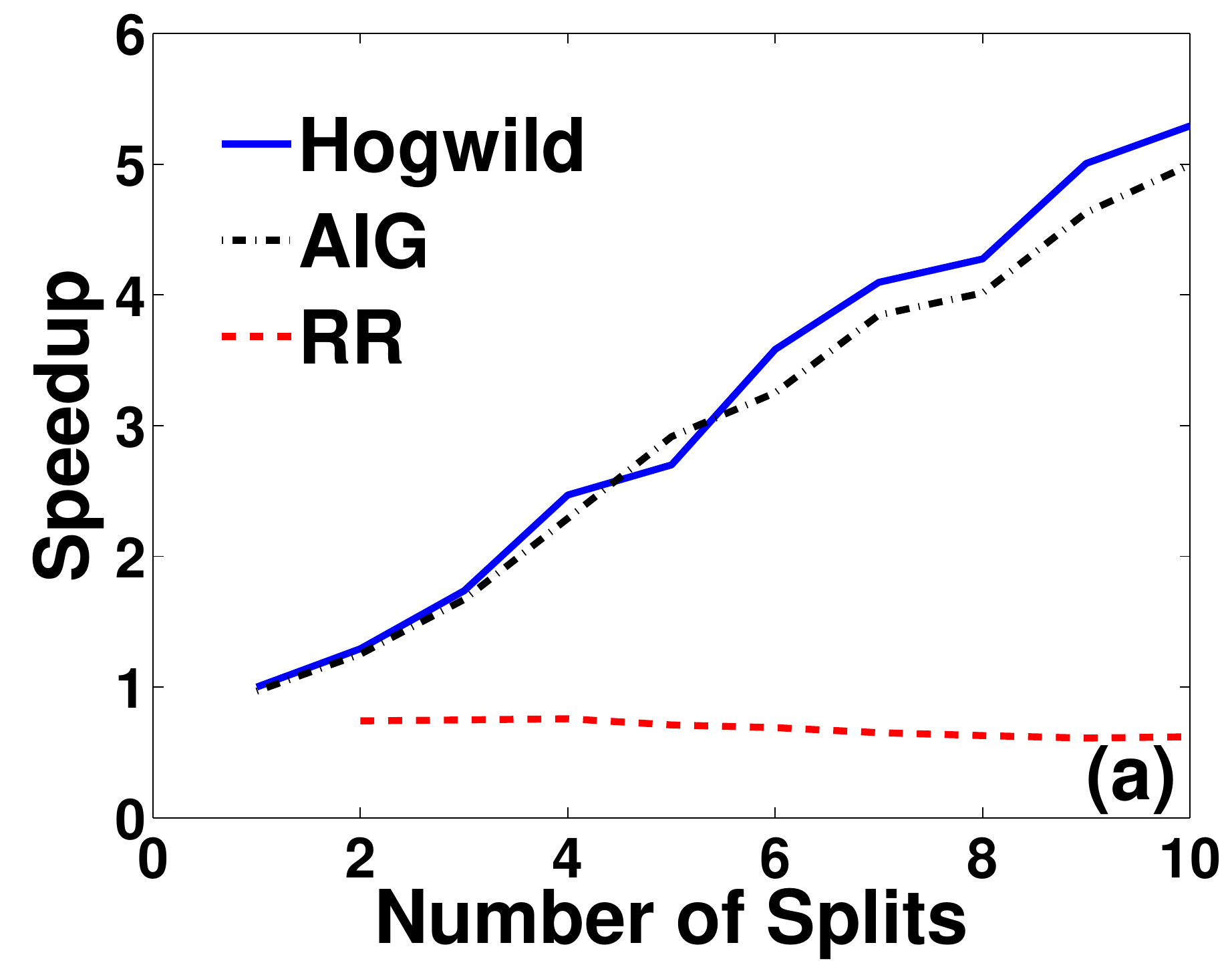} & 
\includegraphics[width=50mm]{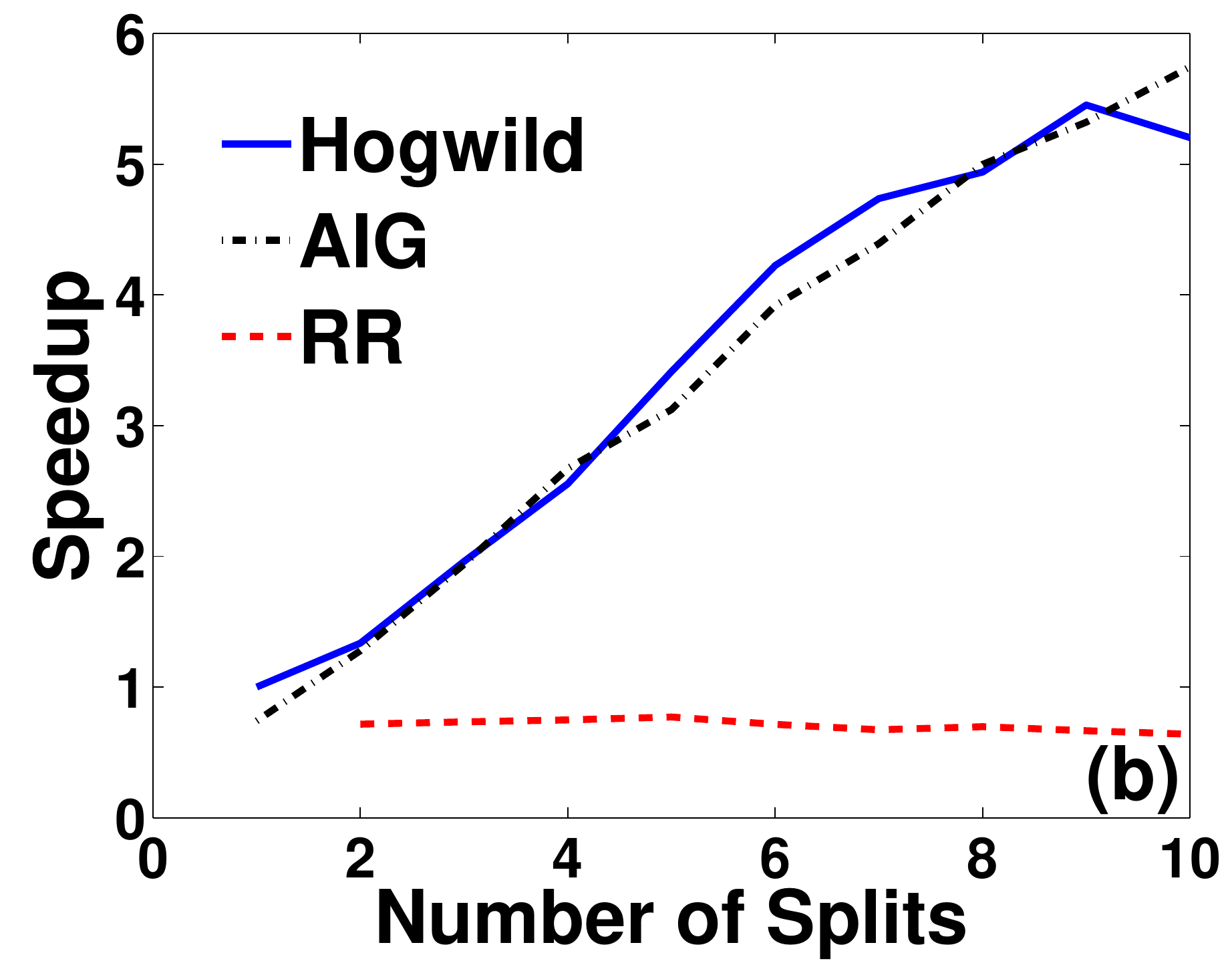} & 
\includegraphics[width=50mm]{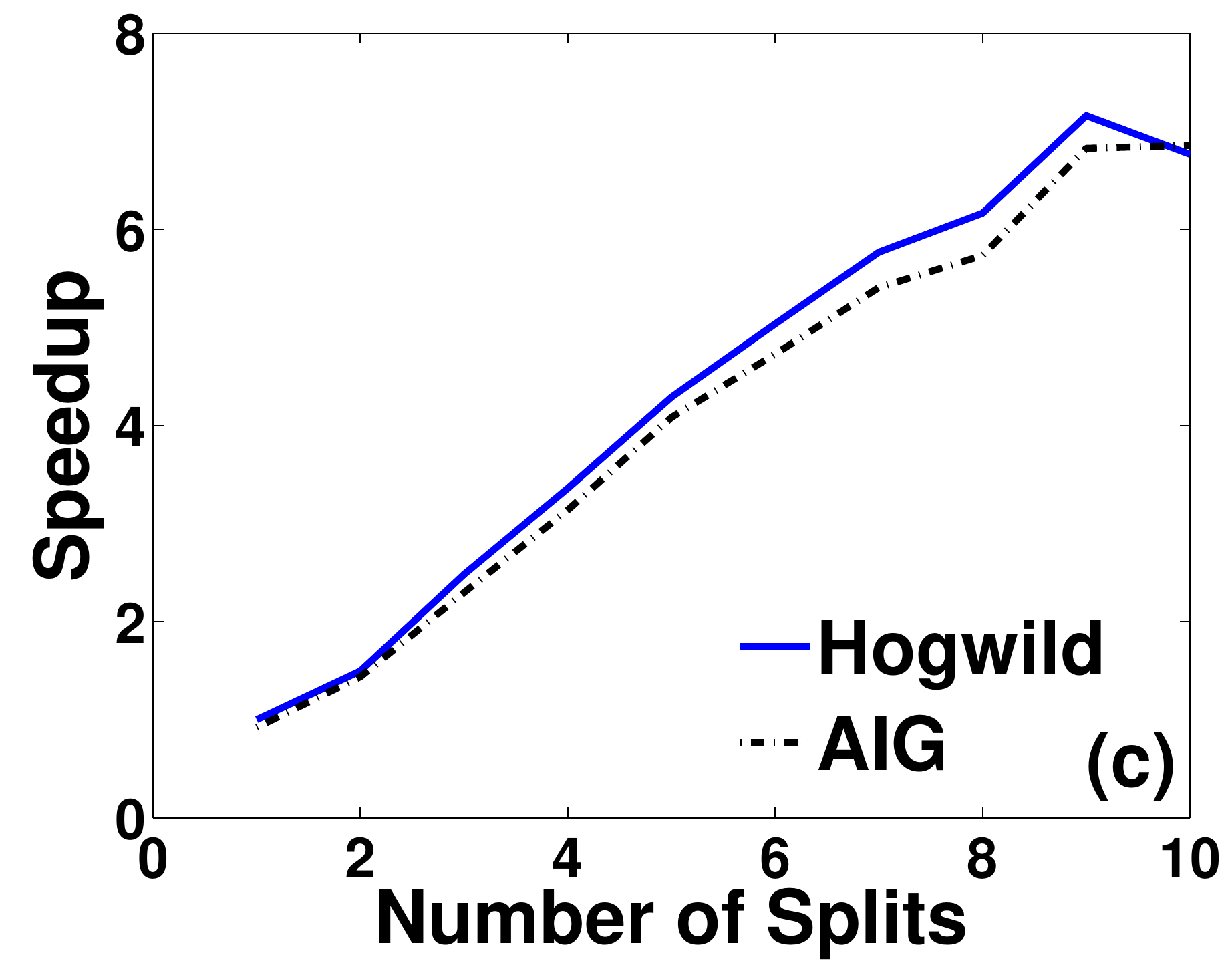} 
\end{tabular}
\caption{Total CPU time versus number of threads for the matrix completion problems (a) Netflix Prize, (b) KDD Cup 2011, and (c) the synthetic Jumbo experiment.
\label{fig:mc-speedups}}
\end{figure}

Figure~\ref{fig:misc}(a) shows the speedups for these three data sets
using \name.  Note that the Jumbo and KDD examples do not fit in our
allotted memory, but even when reading data off disk, \name attains a
near linear speedup. The Jumbo problem takes just over two and a half
hours to complete.  Speedup graphs 
comparing \name to AIG and RR on the three matrix completion
experiments are provided in Figure~\ref{fig:mc-speedups}.  Similar to
the other experiments with quickly computable gradients, RR does not
show any improvement over a serial approach.  In fact, with 10
threads, RR is 12\% slower than serial on KDD Cup and 62\% slower on
Netflix. In fact, it is too slow to complete the Jumbo experiment in any reasonable amount of time, while the 10-way parallel \name implementation solves this problem in under three hours.

\textbf{Graph Cuts.} Our first cut problem was a standard image segmentation by graph cuts problem popular in computer vision.  We computed a two-way cut of the {\tt abdomen} data set~\cite{MaxFlowData}.  This data set consists of a volumetric scan of a human abdomen, and the goal is to segment the image into organs.  The image has $512\times 512\times 551$ voxels, and the associated graph is 6-connected with maximum capacity 10.  Both $\rho$ and $\Delta$ are equal to 9.2e-4  We see that \name speeds up the cut problem by more than a factor of 4 with 10 threads, while RR is twice as slow as the serial version.

Our second graph cut problem sought a mulit-way cut  to determine entity recognition in a large database of web data.  We created a data set of clean entity lists from the DBLife website and of entity mentions from the DBLife Web Crawl~\cite{dblife}.  The data set consists of 18,167 entities and 180,110 mentions and similarities given by string similarity.   In this problem each stochastic gradient step must compute a Euclidean projection onto a simplex of dimension 18,167.  As a result, the individual stochastic gradient steps are quite slow. Nonetheless, the problem is still very sparse with $\rho$=8.6e-3 and $\Delta$=4.2e-3. Consequently, in Figure~\ref{fig:speedups}, we see the that \name achieves a ninefold speedup with 10 cores.  Since the gradients are slow, RR is able to achieve a parallel speedup for this problem, however the speedup with ten processors is only by a factor of 5.   That is, even in this case where the gradient computations are very slow, \name outperforms a round-robin scheme.

\textbf{What if the gradients are slow?} As we saw with the DBLIFE data set, the RR method does get a nearly linear speedup when the gradient computation is slow.  This raises the
question whether RR ever outperforms \name for slow gradients.  To
answer this question, we ran the RCV1 experiment again and introduced
an artificial delay at the end of each gradient computation to
simulate a slow gradient.  In Figure~\ref{fig:misc}(c), we plot the
wall clock time required to solve the SVM problem as we vary the delay
for both the RR and \name approaches.

Notice that \name achieves a greater decrease in computation time
across the board.  The speedups for both methods are the same when the delay is few milliseconds.  That is, if a gradient takes longer than one millisecond to compute,  RR is on par with \name (but not better).  At this rate, one is only able to
compute about a million stochastic gradients per hour, so the gradient computations must be very labor intensive in order for the RR method to be competitive.

\section{Conclusions}

Our proposed \name algorithm takes advantage of sparsity in machine learning problems to enable near linear speedups on a variety of applications.  Empirically, our implementations outperform our theoretical analysis.  For instance, $\rho$ is quite large in the RCV1 SVM problem, yet we still obtain significant speedups.  Moreover,  our algorithms allow parallel speedup even when the gradients are computationally intensive.

Our \name schemes can  be generalized to problems where some of the variables occur quite frequently as well.  We could choose to not update certain variables that would be in particularly high contention.  For instance, we might want to add a bias term to our Support Vector Machine, and we could still run a \name scheme, updating the bias only every thousand iterations or so.

For future work, it would be of interest to enumerate structures that allow for parallel gradient computations with no collisions at all. That is, it may be possible to bias the SGD iterations to completely avoid memory contention between processors.   For example, recent work proposed a biased ordering of the stochastic gradients in matrix completion problems that completely avoids memory contention between processors~\cite{RechtRe11}.   An investigation into how to generalize this approach to other structures and problems would enable even faster computation of machine learning problems.

\section*{Acknowledgements}

BR is generously supported by ONR award N00014-11-1-0723 and NSF award CCF-1139953.  CR is generously supported by the Air Force Research Laboratory (AFRL) under prime contract no. FA8750-09-C-0181, the NSF CAREER award under IIS-1054009, ONR award N000141210041, and gifts or research awards from Google, LogicBlox, and Johnson Controls, Inc.  SJW is generously supported by NSF awards DMS-0914524 and
DMS-0906818 and DOE award DE-SC0002283.  Any opinions, findings, and conclusion or recommendations expressed in this work are those of the authors and do not necessarily reflect the views of any of the above sponsors including DARPA, AFRL, or the US government.

{\small
\bibliographystyle{abbrv}
\bibliography{brecht}
}

\appendix

\section{Analysis of \name}

It follows by rearrangement of~\eq{def.c} that
\begin{equation} \label{eq:hw-app.c1} 
(x-x')^T\nabla f(x) \ge f(x)-f(x') + \frac12 c \|x-x'\|^2, \;\; \mbox{for all  $x \in X$.}
\end{equation}
In particular, by setting $x'=x_\star$ (the minimizer) we have
\begin{equation} \label{eq:hw-app.c2} 
(x-x_\star)^T \nabla f(x)  \ge \frac12 c \|x-x_\star\|^2,
\;\; \mbox{for all $x \in X$.}
\end{equation}
We will make use of these identities frequently in what follows.

\subsection{Developing the recursion}

For $1\leq v \leq n$, let ${\cal P}_v$ denote the projection onto the
$v$th component of $\R^n$. (${\cal P}_v$ is a diagonal matrix with a 1
in the $v$th diagonal position and zeros elsewhere.) Similarly, for
$e\in E$, we let ${\cal P}_e$ denote the projection onto the
components indexed by $e$ --- a diagonal matrix with ones in the
positions indexed by $e$ and zeros elsewhere. 

We start with the update formula \eq{hw-wr-update}.  Recall that $k(j)$ is the state of the decision variable's counter when the update to $x_j$ was read. We have
\[
	x_{j+1} = x_j
-  \gamma |e_j| {\cal P}_{v_j} G_{e_j}(x_{k(j)})\,.
\]

By subtracting $x_\star$
from both sides, taking norms, we have
\begin{align*}
\frac12 \|x_{j+1} - x_\star\|_2^2  &= \frac12 \|x_{j}-x_\star\|_2^2 
- \gamma |e_j| (x_j-x_\star)^T  {\cal P}_{v_j}G_{e_j}(x_{k(j)}) \\
&\qquad\qquad + \frac12 \gamma^2 |e_j|^2 \|{\cal P}_{v_j}G_{e_j}(x_{k(j)})\|^2 \\
& =  \frac12 \|x_{j}-x_\star\|_2^2 - \gamma |e_j|  (x_j-x_{k(j)})^T  {\cal P}_{v_j}G_{e_j}(x_{j}) \\
&\qquad  -\gamma |e_j|  (x_j-x_{k(j)})^T {\cal P}_{v_j} (G_{e_j}(x_{k(j)})- G_{e_j}(x_{j})) \\
&\qquad\qquad- \gamma |e_j|  (x_{k(j)}-x_\star)^T  {\cal P}_{v_j} G_{e_j}(x_{k(j)})  +\frac12 \gamma^2 |e_j|^2 \|{\cal P}_{v_j}G_{e_j}(x_{k(j)})\|^2 
\end{align*}
Let $a_j = \frac12 \E[\|x_j - x_\star\|^2_2]$. By taking expectations of both sides and using the bound (\ref{eq:def.M}), we obtain
\begin{equation} \label{eq:hw-Ajchange}
\begin{aligned}
a_{j+1} &\le a_j - \gamma \E[(x_j-x_{k(j)})^T G_{e_j}(x_{j})] 
- \gamma \E[(x_j-x_{k(j)})^T (G_{e_j}(x_{k(j)})-G_{e_j}(x_{j}))] \\
&\qquad \qquad - \gamma \E[(x_{k(j)}-x_\star)^T G_{e_j}(x_{k(j)})]  +\frac12 \gamma^2 \Omega M^2.
\end{aligned}
\end{equation}
where we recall that $\Omega = \max_{e\in E} |e|$.  Here, in several places, we used the useful identity: for any function $\zeta$ of $x_1,\ldots, x_j$ and any $i \leq j$, we have
\[
\begin{aligned}
	\E[|e_j| {\cal P}_{v_j} G_{e_j}(x_i)^T \zeta(x_1,\ldots,x_j)] &= 
		\E\left\{ \E\left[|e_j| {\cal P}_{v_j} G_{e_j}(x_i)^T \zeta(x_1,\ldots,x_j)~|~e_1,\ldots,e_j,v_1,\ldots, v_{j-1}\right] \right\}\\
		& = \E\left[ G_{e_j}(x_i)^T \zeta(x_1,\ldots,x_j)\right] \,.
		\end{aligned}
\]
We will perform many similar calculations throughout, so, 
before
proceeding, we denote
\[
e_{[i]} := (e_1,e_2,\dotsc,e_i, v_1,v_2,\dotsc,v_i),
\]
to be the tuple of all edges and vertices selected in updates $1$ through $i$. Note that $x_l$ depends on $e_{[l-1]}$ but not on $e_j$  or $v_j$ for any $j \ge l$.   We next consider the three expectation terms in this expression.

Let's first bound the third expectation term in \eq{hw-Ajchange}.   Since $x_{k(j)}$ is independent of $e_j$ we have 
\begin{align*}
\E & \left[ (x_{k(j)}-x_\star)^T G_{e_j}(x_{k(j)}) \right]  \\
&= \E \left\{ \E \left[ (x_{k(j)}-x_\star)^T G_{e_j}(x_{k(j)}) 
| e_{[k(j)-1]}\right] \right\} \\
&= \E \left\{ (x_{k(j)}-x_\star)^T \E \left[  G_{e_j}(x_{k(j)}) 
| e_{[k(j)-1]}\right] \right\} \\ 
&= \E \left[ (x_{k(j)}-x_\star)^T \nabla f (x_{k(j)}) \right]\,,
\end{align*}
where $\nabla f$ denotes an element of $\partial f$.  It follows from
(\ref{eq:hw-app.c2}) that
\begin{equation} \label{eq:hw-Efirst}
\E \left[ (x_{k(j)}-x_\star)^T \nabla f(x_{k(j)}) \right] \ge 
c a_{k(j)}\,.
\end{equation}

The first expectation can be treated similarly:
\begin{align}
\nonumber
	\E[(x_j-x_{k(j)})^T G_{e_j}(x_{j})]  
	&= \E\left\{ \E[(x_j-x_{k(j)})^T G_{e_j}(x_{j})~|~e_{[j-1]}]\right\}\\
\nonumber
	&= \E\left\{(x_j-x_{k(j)})^T  \E[G_{e_j}(x_{j})~|~e_{[j-1]}]\right\}\\
\nonumber
	&=  \E[(x_j-x_{k(j)})^T \nabla  f(x_{j})] \\
\label{eq:jkj1}
	&\ge \E[f(x_j)-f(x_{k(j)})]+ \frac{c}{2} \E[\|x_j-x_{k(j)}\|^2]
\end{align}
where the final inequality is from~\eq{hw-app.c1}.  Moreover, we can estimate the difference between $f(x_j)$ and $f(x_{k(j)})$ as 
\begin{align}
\nonumber
	\E[ f(x_{k(j)}) - f(x_j) ] &= \sum_{i=k(j)}^{j-1} \E[ f(x_{i}) - f(x_{i+1})] \\
\nonumber
	& = \sum_{i=k(j)}^{j-1} \sum_{e\in E}  \E[ f_e(x_{i}) - f_e(x_{i+1})] \\
\nonumber
	& \leq \frac{\gamma}{|E|} \sum_{i=k(j)}^{j-1} \sum_{e\in E}  \E [ G_e(x_i)^T G_{e_i}(x_i)]  \\
\label{eq:jkj2}
	& \leq \gamma \tau   \rho M^2 \,.
\end{align}
Here we use the inequality
\[
	f_e(x_{i}) - f_e(x_{i+1}) \leq \frac{1}{|E|} G_e(x_{i})^T(x_i - x_{i+1}) = \frac{\gamma}{|E|} G_e(x_i)^T G_{e_i}(x_i)
\]
which follows because $f_e$ is convex.  By combining (\ref{eq:jkj1}) and
(\ref{eq:jkj2}), we obtain
\begin{equation}
\label{eq:jkj3}
	\E[(x_j-x_{k(j)})^T G_{e_j}(x_{j})] \ge -\gamma \tau   \rho M^2+ \frac{c}{2} \E[\|x_j-x_{k(j)}\|^2]\,.
\end{equation}

We turn now to the second expectation term in \eq{hw-Ajchange}.  We have
\begin{alignat}{2}
\nonumber
 \E &  \left[ (x_j-x_{k(j)})^T (G_{e_j}(x_{k(j)})-G_{e_j}(x_j)) \right]  & & \\
\nonumber
&= \E \left[ \sum_{i=k(j)}^{j-1}  (x_{i+1}-x_i)^T (G_{e_j}(x_{k(j)})-G_{e_j}(x_j))  \right]   && \\
\nonumber
&=  \E \left[  \sum_{i=k(j)}^{j-1} \gamma |e_i| G_{e_i}(x_{k(i)})^T (G_{e_j}(x_{k(j)})-G_{e_j}(x_j))  \right]   \\
\nonumber
&=  \E \left[  \sum_{\stackrel{i=k(j)}{e_i \cap e_j \neq \emptyset}}^{j-1} \gamma |e_i| G_{e_i}(x_{k(i)})^T (G_{e_j}(x_{k(j)})-G_{e_j}(x_j)) \right]   & & \\
\nonumber
&\ge -\E \left[ \sum_{\stackrel{i=k(j)}{e_i \cap e_j \neq \emptyset}}^{j-1} \gamma |e_i| \| G_{e_i}(x_{k(i)}) \| \, \| G_{e_j}(x_{k(j)})-G_{e_j}(x_j) \| \right] && \\
\nonumber
& \ge -\E \left[ \sum_{\stackrel{i=k(j)}{e_i \cap e_j \neq \emptyset}}^{j-1} 2\Omega M^2 \gamma \right] \\
\label{eq:jkj4}
& \ge -2\Omega M^2  \gamma \rho\tau
\end{alignat}
where $\rho$ is defined by~\eq{collide-prob}.  Here, the third line follows from our definition of the gradient update.  The fourth line is tautological: only the edges where $e_i$ and $e_j$ intersect nontrivially factor into the sum.  The subsequent inequality is Cauchy-Schwarz, and the following line follows from~\eq{def.M}.

By substituting \eq{hw-Efirst}, (\ref{eq:jkj3}), and (\ref{eq:jkj4}) into
\eq{hw-Ajchange}, we obtain the following
bound:
\begin{align}
\label{eq:hw-Ajred.0}
a_{j+1} & \le a_j -  c\gamma \left(a_{k(j)} + \tfrac{1}{2} \E[\|x_j-x_{k(j)}\|^2]\right) 
+ \frac{M^2 \gamma^2}{2}  \left( \Omega + 2\tau \rho + 4 \Omega \rho \tau \right)\,.
\end{align}

To complete the argument, we need to bound the remaining expectation in~\eq{hw-Ajred.0}.  We expand out the expression multiplied by $c\gamma$ in~\eq{hw-Ajred.0} to find
\begin{align*}
	a_{k(j)} +\tfrac{1}{2} \E[\|x_j-x_{k(j)}\|^2 &= a_j - \E\left[(x_{j}-x_{k(j)})^T (x_{k(j)}-x_\star)\right]\\
	&= a_j  - \E\left[ \sum_{i=k(j)}^{j-1} (x_{i+1}-x_{i})^T (x_{k(j)}-x_\star) \right]\\
	&= a_j  - \E\left[ \sum_{i=k(j)}^{j-1} \gamma |e_i|G_{e_i}(x_{k(i)})^T {\cal P}_{v_i}(x_{k(j)}-x_\star) \right]\,.
\end{align*}

Let $e_{[\neg i ]}$ denote the set of all sampled edges and vertices except for
$e_i$ and $v_i$.  Since $e_i$ and $v_i$ are both independent of $x_{k(j)}$, we can proceed to
bound
\begin{align*}
&\E\left[ \sum_{i=k(j)}^{j-1} \gamma |e_i| G_{e_i}(x_{k(i)})^T{\cal P}_{v_i} (x_{k(j)}-x_\star) \right]\\
	\leq&  \E\left[ \sum_{i=k(j)}^{j-1} \gamma \Omega M \|{\cal P}_{v_i} (x_{k(j)}-x_\star)\|_2 \right]\\
	=& \gamma \Omega M  \sum_{i=k(j)}^{j-1} \E\left[ \E_{e_i,v_i}\left[\|{\cal P}_{e_i,v_i} (x_{k(j)}-x_\star)\|_2~|~e_{[\neg i]} \right] \right]\\	
	\leq&   \gamma \Omega M  \sum_{i=k(j)}^{j-1} \E\left[ \left(\E_{e_i,v_i}\left[ (x_{k(j)}-x_\star)^T {\cal P}_{v_i} (x_{k(j)}-x_\star)~|~e_{[\neg i]}\right]\right)^{1/2}\right]\\
	= & \gamma \Omega M  \sum_{i=k(j)}^{j-1} \E\left[\left((x_{k(j)}-x_\star)^T \E_{e,v}[{\cal P}_{v}] (x_{k(j)}-x_\star)\right)^{1/2}\right]\\	
	\leq&  \tau \gamma\Omega M  \E\left[\left((x_{k(j)}-x_\star)^T \E_{e,v}[{\cal P}_{v}] (x_{k(j)}-x_\star)\right)^{1/2}\right]\\
	\leq&  \tau \gamma \Omega M  \Delta^{1/2}  \E[\| x_{k(j)}-x_\star\|_2] \\
	\leq&  \tau \gamma \Omega M \Delta^{1/2} \left(\E[\| x_{j}-x_\star\|_2]  + \tau \gamma \Omega M\right)\\
	\leq& \tau \gamma \Omega M \Delta^{1/2} (\sqrt{2}a_j^{1/2} + \tau \gamma \Omega M)\, , 
\end{align*}
where $\Delta$ is defined in~\eq{collide-prob}. The first inequality is Cauchy-Schwartz.  The next inequality is Jensen.  The second to last inequality follows from our definition of $x_j$, and the final inequality is Jensen again.

Plugging the last two expressions into~\eq{hw-Ajred.0}, we obtain
\begin{equation}\label{eq:hw-final-recurse}
a_{j+1}  \le (1 -  c\gamma) a_j 
+ \gamma^2 \left(\sqrt{2} c \Omega M\tau \Delta^{1/2}\right) a_j^{1/2} +\frac12
 M^2 \gamma^2  Q
\end{equation}
where 
\[
	Q= \Omega + 2\tau \rho + 4 \Omega \rho \tau + 2\tau^2 \Omega^2 \Delta^{1/2}\,.
\]
Here we use the fact that $c\gamma<1$ to get a simplified form for $Q$.  This recursion only involves constants involved with the structure of $f$, and the nonnegative sequence $a_j$.  To complete the analysis, we will perform a linearization to put this recursion in a more manageable form.

To find the steady state, we must solve the equation
\begin{equation}\label{eq:hw-ss}
a_\infty  = (1 -  c\gamma) a_\infty
+ \gamma^2 \left(\sqrt{2} c \Omega M\tau \Delta^{1/2}\right) a_\infty^{1/2} 
+  \frac{M^2\gamma^2}{2} Q\,.
\end{equation}
This yields the fixed point
\begin{equation}\label{eq:fixed-point}
\begin{aligned}
	a_\infty &=  \frac{M^2\gamma^2}{2}\left( \Omega \tau \Delta^{1/2} + \sqrt{\Omega^2 \tau^2 \Delta^ + \frac{Q}{c\gamma}}\right)^2 \\
	&\leq \frac{M^2\gamma}{2c}\left( \Omega \tau \Delta^{1/2} + \sqrt{\Omega^2 \tau^2 \Delta + Q}\right)^2 = C(\tau,\rho,\Delta,\Omega) \frac{M^2\gamma}{2c}	\,.
	\end{aligned}
\end{equation}
Note that for $\rho$ and $\Delta$ sufficiently small, $C(\tau,\rho,\Delta,\Omega)\approx 1$.

Since the square root is concave, we can linearize~\eq{hw-final-recurse} about the fixed point $a_\infty$ to yield
\[
\begin{aligned}
a_{j+1} & \le (1 -  c\gamma) (a_j-a_\infty) 
+  \frac{\gamma^2 \left(\sqrt{2} c \Omega M\tau \Delta^{1/2}\right)}{2\sqrt{a_\infty}} {(a_j - a_\infty)}\\
&\qquad  +  (1-c\gamma)a_\infty
+ \gamma^2 \left(\sqrt{2} c \Omega M\tau \Delta^{1/2}\right) a_\infty^{1/2} 
+  \frac12 M^2\gamma^2 Q\\
&=   \left(1 -  c\gamma ( 1    -  \delta ) \right)(a_j-a_\infty) 
+  a_\infty\,.
\end{aligned}
\]
Here
\[
	\delta =  \frac{1}{1 + \sqrt{ 1 + \frac{Q}{c\gamma \Omega^2 \tau^2 \Delta}}} \leq \frac{1}{1 + \sqrt{ 1 + \frac{Q}{\Omega^2 \tau^2 \Delta}}} 
\]

To summarize, we have shown that the sequence $a_j$ of squared distances satisfies
\begin{equation}\label{eq:hw-linearized-final-recurse}
	a_{j+1} \leq (1-c\gamma(1-\delta(\tau,\rho,\Delta,\Omega))) (a_j - a_\infty) + a_\infty	
\end{equation}
with $a_\infty \leq C(\tau,\rho,\Delta,\Omega) \frac{M^2\gamma}{2 c}$.    In the case that $\tau = 0$ (the serial case), $C(\tau,\rho,\Delta,\Omega)=\Omega$ and $\delta(\tau,\rho,\Delta,\Omega)=0$.  Note that if $\tau$ is non-zero, but $\rho$ and $\Delta$ are $o(1/n)$ and $o(1/\sqrt{n})$ respectively, then as long as $\tau = o(n^{1/4})$, $C(\tau,\rho,\Delta,\Omega)=O(1)$.  In our setting, $\tau$ is proportional to the number of processors, and hence as long as the number of processors is less $n^{1/4}$, we get nearly the same recursion as in the linear rate.

In the next section, we show that~\eq{hw-linearized-final-recurse} is sufficient to yield a $1/k$ convergence rate.  Since we can run $p$ times faster in our parallel setting, we get a linear speedup.

\subsection{Proof of Proposition~\ref{prop:hw-main-result}: Final Steps}

Since $\nabla f$ is Lipschitz, we have
\[
	f(x)\leq f(x') + \nabla f(x')^T(x-x') + \frac{L}{2}\|x-x'\|^2\quad \mbox{for all} ~x, x'\in\R^n\,.
\] 
Setting $x'=x_\star$ gives $f(x)-f(x_\star) \leq \frac{L}{2} \|x-x_\star\|^2$.  Hence,
\[
	\E[f(x_k) - f(x_\star)] \leq L a_k
\]
for all $k$.  To ensure the left hand side is less than $\epsilon$, it suffices to guarantee that $a_k \leq \epsilon/L$.

To complete the proof of Proposition~\ref{prop:hw-main-result}, we use the results of Section~\ref{sec:fast-rates-theory}.  We wish to achieve a target accuracy of $\epsilon/L$.  To apply~\eq{hw-one-epoch}, choose $a_\infty$ as in~\eq{fixed-point} and the values
\[
\begin{aligned}
	c_r &= c(1-\delta(\tau,\rho,\Delta,\Omega)\\
	 B &=  C(\tau,\rho,\Delta,\Omega) \frac{M^2\gamma}{2c}\,.
\end{aligned}
\]
By~\eq{fixed-point}, we have $a_\infty \leq \gamma B$.

Choose $\gamma$ satisfying~\eq{hw-gamma-ub}.  With this choice, we automatically have
\[
	\gamma \leq  \frac{\epsilon c}{LM^2\Omega^2 \tau^2\Delta \left(1 + \sqrt{1+\frac{Q}{\Omega^2\tau^2\Delta}}\right)^2} \leq \frac{\epsilon}{2LB}
\]
because $(1+\sqrt{1+x})^2 \leq 4+2x$ for all $x\geq 0$. Substituting this value of $\gamma$ into~\eq{hw-one-epoch}, we see that
\begin{equation}\label{eq:hw-final-iter-ub}
	k \geq \frac{LM^2\log(LD_0/\epsilon)}{\epsilon c^2} \cdot\frac{C(\tau,\rho,\Delta,\Omega)}{1-\delta(\tau,\rho,\Delta,\Omega)}
\end{equation}
iterations suffice to achieve $a_k \leq \epsilon/L$.  Now observe that
\[
\begin{aligned}
\frac{C(\tau,\rho,\Delta,\Omega)}{1-\delta(\tau,\rho,\Delta,\Omega)} &=
\Omega^2 \tau^2\Delta \cdot \frac{\left(1+\sqrt{1+\frac{ Q }{\Omega^2 \tau^2\Delta}}\right)^2}{1 - \frac{1}{1+\sqrt{1+\frac{ Q}{\Omega^2 \tau^2\Delta}}}}\\
& = \Omega^2 \tau^2 \Delta\cdot \frac{\left(1+\sqrt{1+\frac{ Q}{\Omega^2 \tau^2\Delta}}\right)^2}{\sqrt{1+\frac{Q}{\Omega^2\tau^2\Delta}}}\\
&\leq  8 \Omega^2 \tau^2 \Delta +2Q\\
&=  8 \Omega^2 \tau^2 \Delta + 2\Omega + 4\tau \rho + 8 \Omega \rho \tau + 4\tau^2 \Omega^2 \Delta^{1/2}\\
&\leq  2\Omega ( 1+ 6\tau \rho + 6 \tau^2 \Omega \Delta^{1/2})   \,.
\end{aligned}
\]
Here,the second to last inequality follows because 
\[
	\frac{\left(1+\sqrt{1+x}\right)^3}{\sqrt{1+x}} \leq 8 + 2x
\]
for all $x\geq 0$.  Plugging this bound into~\eq{hw-final-iter-ub} completes the proof.

\end{document}